\documentclass[11pt,a4paper]{article}
\usepackage{latexsym,amssymb,amsmath,amsthm,amsfonts,enumerate,verbatim,xspace,
exscale}

\input xy
\xyoption{all}
\CompileMatrices
\UseComputerModernTips

\usepackage{graphicx}
\usepackage{amsmath,amsbsy,amsfonts,amssymb}

\title{Reduction of Almost Poisson brackets and Hamiltonization of the Chaplygin Sphere}
\author{{\sc{Luis C. Garc\'ia-Naranjo} \thanks{
          Section de Math\'ematiques, Ecole Polytechnique F\'ed\'erale de Lausanne, CH-1015 
Lausanne, Switzerland. \newline{\texttt{E-mail: luis.garcianaranjo@epfl.ch}}}}}

\theoremstyle{plain}
\newtheorem{theorem}{Theorem}[section]

\newtheorem{proposition}[theorem]{Proposition}
\newtheorem{corollary}[theorem]{Corollary}

\theoremstyle{definition}
\newtheorem*{definition}{Definition}

\newenvironment{remark}[1][Remark]{\begin{trivlist}
\item[\hskip \labelsep {\bfseries #1}]}{\end{trivlist}}

\newcommand{\g}{\mathfrak{g}}

\def\so{\mathfrak{so}}
\def\se{\mathfrak{se}}

\def\J{\mathbf{J}}

\def\g{\mathfrak{g}}
\def\h{\mathfrak{h}}
\def\D{\mathcal{D}}
\def\R{\mathbb{R}}
\def\I{\mathbb{I}}
\def\L{\mbox{Leg}}
\def\M{\mathcal{M}}

\def\F{\mathcal{F}}

\def\Ham{\mathcal{H}}

\def\Lag{\mathcal{L}}

\def\C{\mathcal{C}}

\newcommand{\up}{\upshape}

\def\vv<#1>{\langle#1\rangle}

\providecommand{\det}{\mbox{$\text{\up{det}}\,$}}

\newcommand{\Ad}{\mbox{$\text{\upshape{Ad}}$}}

\begin{document}
\maketitle

\begin{abstract}
We construct different  almost Poisson brackets for nonholonomic systems than those existing in the literature and study their reduction.    
Such  brackets are built by considering   non-canonical two-forms  on the cotangent bundle of configuration space and then 
carrying out a  projection  onto the constraint space that encodes the Lagrange-D'Alembert principle.
We justify the need for this type of brackets  by working out the reduction of
the celebrated Chaplygin sphere rolling problem.
 Our construction provides a
geometric explanation of the Hamiltonization of the
problem given by A. V. Borisov and I. S. Mamaev.
\end{abstract}


\section{Introduction and Outline}

The equations of motion for nonholonomic systems are not Hamiltonian. They can be formulated 
with respect to an almost nonholonomic Poisson bracket of functions that fails to satisfy the Jacobi identity. This 
formulation has its origins in \cite{vanderschaft1994, Marle1998, Cantrijn99}, and others. Roughly speaking, one constructs the  bracket by projecting the canonical Poisson tensor on the cotangent bundle of configuration space onto the constraint space using the Lagrange-D'Alembert principle.  This formulation is of interest because the reduced equations of motion of some important  examples have
been written in Hamiltonian form with respect a usual Poisson bracket (sometimes after a time rescaling), \cite{{BorisovMamaev},
{BorisovMamaev_Rolling_Body1}, {BorisovMamaev_Rolling_Body2}, {BorisovMamaev1995},{Fasso},{Naranjo2007}}.  In this case we  say that the problem has been \emph{Hamiltonized}. An important example of Hamiltonization, although independent of bracket formulations is given in \cite{FedorovJovan}.

In this 
paper we construct  more general almost Poisson structures and study their reduction.
These
brackets are obtained by projecting \emph{non-canonical} bi-vector fields on the cotangent bundle of configuration
space onto the constraint space using the Lagrange-D'Alembert principle. 
We call these brackets  \emph{Affine Almost Poisson Brackets} since they are derived from a 
non-degenerate two-form  that deviates
from the canonical one by the addition of an \emph{affine} term of magnetic type that annihilates the free Hamiltonian vector field. The idea of adding an 
affine term that ``does not see the flow" already appears in \cite{EhlersKoiller}.

Our motivation for considering this type of brackets was to link the general construction 
of almost Poisson brackets, \cite{vanderschaft1994, Marle1998, Cantrijn99}, with the Hamiltonian formulation of the equations of motion for the celebrated Chaplygin sphere problem given in \cite{BorisovMamaev, BorisovMamaev1995}. Although the authors in \cite{EhlersKoiller}  consider affine symplectic structures, the 
Hamiltonization of the problem could not be verified by their methods.

Perhaps the most important result in the context of Hamiltonization is Chaplygin's reducing
multiplier theory \cite{Chapligyn_reducing_multiplier} that applies to nonholonomic systems with two degrees of freedom that posses an invariant measure. This theory does not apply directly to the Chaplygin sphere problem since the reduced space is 5-dimensional. However, associated to the symmetry there is a conserved quantity, the vertical angular momentum. This integral is used by the authors in \cite{BorisovMamaev1995} to apply Routh's method of reduction to the reduced equations and then apply a generalized version of Chaplygin's theory to Hamiltonize the problem.

We give an alternative method to achieve the Hamiltonization of the Chaplygin sphere by carrying out the reduction of an affine almost Poisson bracket. The affine term in the bracket needs to be included for the conserved quantity to become a Casimir function of the corresponding reduced bracket. In fact, we also show that the reduction of the standard nonholonomic bracket defined in  \cite{vanderschaft1994, Marle1998} does not even yield a foliation of the reduced space by even dimensional leaves, having thus
very different properties from a usual Poisson bracket. 
This shows that one should extend the notion of  nonholonomic almost Poisson brackets as introduced 
in \cite{vanderschaft1994, Marle1998}, and generally considered in the literature, and incorporate  the affine description if one is interested in obtaining reduced brackets with optimal properties.

The outline of the paper is as follows. After reviewing the construction of almost Poisson brackets in section \ref{S:Background}, we present a basic scheme for their reduction in section \ref{S:Reduction_of_Almost_Poisson_Brackets}. In section \ref{S:Affine_brackets} we introduce the notion
of affine almost Poisson brackets and discuss their reduction in section \ref{S:Reduction_of_Affine_Nonholonomic_Brackets}. Section \ref{S:Chaplygin_Sphere}
treats the Chaplygin sphere problem. We use Cartan's moving
frames for $SO(3)$ to avoid  messy calculations in Euler angles or other local coordinates.
Finally, we give some closing remarks in section \ref{S:Final_Remarks}.
\section{Background}
\label{S:Background}

\subsection*{Nonholonomic Systems}
\label{S:Nonholonomic_Systems}

A nonholonomic system consists of  a configuration space $Q$ with local coordinates $q_a$, $a=1,\dots, n$, a hyper-regular Lagrangian $\Lag:TQ\rightarrow \R$, and 
 a non-integrable distribution $\D\subset TQ$ that describes
the kinematic nonholonomic constraints. 
In coordinates the distribution is defined by the independent
equations\footnote{Here and in what follows the Einstein convention of sum over
repeated indices holds and we assume smoothness of all quantities.}
\begin{equation}
\label{E:Constraints_Lag_form}
\epsilon^i_a(q)\dot q_a=0, \qquad i=1,\dots,k<n,
\end{equation}
where the functions $\epsilon^i_a(q)$ are the components of 
the independent \emph{constraint one-forms} on $Q$,
$
\epsilon^i:=\epsilon^i_a(q)\; dq_a.
$

The dynamics of the
system are governed by the
Lagrange-D'Alembert principle. 
This principle states that the forces of constraint annihilate any virtual
displacement so they perform no work during the motion. The equations of motion
take the form 
\begin{equation}
\label{E:Eqns_of_motion_Lag_form}
\frac{d}{dt} \left (\frac{\partial \Lag}{\partial \dot q_a}\right )
-\frac{\partial \Lag}{\partial q_a}=\lambda_i\epsilon^i_a, \qquad a=1,\dots,n.
\end{equation}
The scalar functions $\lambda_i, i=1,\dots,k,$ are referred to as Lagrange
multipliers, that under the assumption of hyper-regularity of the Lagrangian, 
are uniquely determined by the condition that the
constraints (\ref{E:Constraints_Lag_form}) are satisfied.

The equations (\ref{E:Eqns_of_motion_Lag_form}) together with the 
constraints (\ref{E:Constraints_Lag_form})  define a vector field 
$\mbox{$Y_{\textup{nh}}^{\D}$}$ on $\D$  whose integral curves describe the motion of the nonholonomic
system. 
A short calculation shows that along the flow of $\mbox{$Y_{\textup{nh}}^{\D}$}$, the energy function 
$E_\Lag:= \frac{\partial \Lag}{\partial \dot q_a} \; \dot q_a - \Lag$, 
is conserved.

We now write the equations of motion as first order equations on the cotangent bundle $T^*Q$. 
Via the Legendre transform, 
$\L:TQ\rightarrow T^*Q$, we define canonical
coordinates $(q_a,p_a)$ on $T^*Q$ by the rule $\L:(q_a,\dot q_a)\mapsto
(q_a,p_a=\partial \Lag /\partial q_a)$. The Legendre transform is a global diffeomorphism by our
assumption that $\Lag$ is hyper-regular.

The Hamiltonian function, $\Ham:T^*Q\rightarrow \R$, is defined  in the 
usual way $\Ham:=E_\Lag \circ \L^{-1}$. The 
equations of motion (\ref{E:Eqns_of_motion_Lag_form}) are shown to be equivalent
to 
\begin{equation}
\label{E:Eqns_of_motion_Ham_form}
\dot q_a = \frac{\partial \Ham}{\partial p_a}, \qquad 
\dot p_a= - \frac{\partial \Ham}{\partial q_a} + \lambda_i \epsilon_a^i(q), \qquad
a=1,\dots,n,
\end{equation}
and the constraint equations (\ref{E:Constraints_Lag_form}) become
\begin{equation}
\label{E:Constraints_Ham_form}
\epsilon^i_a(q) \frac{\partial \Ham}{\partial p_a}=0, \qquad i=1,\dots,k.
\end{equation}
The above equations define the \emph{constraint submanifold} $\M=\L(\D) \subset
T^*Q$.  Since the Legendre transform is linear on the fibers, $\M$ is a vector
sub-bundle of $T^*Q$ that for each 
$q\in Q$ specifies an $n-k$ vector subspace of $T^*_qQ$.

Equations (\ref{E:Eqns_of_motion_Ham_form}) together with
(\ref{E:Constraints_Ham_form})
define the vector field
$\mbox{$X_{\textup{nh}}^{\M}$}$ on $\M$, that  describes the motion of our nonholonomic system in the
Hamiltonian side and is the push forward of the vector field $\mbox{$Y_{\textup{nh}}^{\D}$}$ by the
Legendre transform. The equations  (\ref{E:Eqns_of_motion_Ham_form}) can be
intrinsically written as
\begin{equation}
\label{E:Eqns_of_motion_Ham_form_intrinsic}
{\bf i}_{\mbox{$X_{\textup{nh}}^{\M}$}} \,\iota^* \Omega_Q = \iota^*(d\Ham +\lambda_i  \tau^* \epsilon^i),
\end{equation}
where $\Omega_Q$ is the canonical symplectic form on $T^*Q$, $\iota:\M
\hookrightarrow T^*Q$ is the inclusion and $\tau: T^*Q\rightarrow Q$ is the
canonical projection. The constraints (\ref{E:Constraints_Ham_form}) and their derivatives are
intrinsically written as
\begin{equation}
\mbox{$X_{\textup{nh}}^{\M}$}\in \C :=T\M\cap \F,
\end{equation}
where $\F$ is the distribution on $T^*Q$ defined as $\F:=\{ v\in T(T^*Q):
\langle \tau^*\epsilon^i,  v\rangle =0\}$.
Here $\C$ is a non-integrable distribution on $\M$  that is denoted by $H$
in \cite{BS93} and $\Ham$ in \cite{KoonMarsden1}.
We reserve these symbols for other objects.

The following theorem was first stated in \cite{Weber1986} and its proof can be found  in \cite{BS93}. 

\begin{theorem}
\label{T:Bates_Snyaticki}
The point-wise restriction of $\iota^* \Omega_Q$ to $\C$, denoted $\Omega_\C$, is non-degenerate.
\end{theorem}

Equivalently we can say that along $\M$ we have the decomposition of $T(T^*Q)$ as 
$T_\M(T^*Q)=\C\oplus \C^{\Omega_Q}$, where $\C^{\Omega_Q}$ is the symplectic orthogonal complement of $\C$ with respect to $\Omega_Q$.

\subsection*{The Almost Hamiltonian Approach}

Although energy is  conserved,  due to the nonholonomic constraints, the equations of motion (\ref{E:Eqns_of_motion_Lag_form}) cannot be cast in Hamiltonian form.
One can however write them with respect to a bracket of functions that fails to
satisfy the Jacobi identity, a so-called \emph{almost Poisson bracket}, see for
example \cite{vanderschaft1994, Marle1998, Cantrijn99}. There is
also an almost symplectic counterpart. After factorization of external
symmetries of a so-called $G$-Chaplygin system,  the equations of motion can be written
with respect to an \emph{almost symplectic form} which is a
non-degenerate two-form that is not closed, see \cite{BS93, KoonMarsden1, Koiller2002}.

Roughly speaking, the process of writing the equations of motion for a
non-holonomic system 
in an almost Hamiltonian way amounts to eliminating the Lagrange multipliers
from the equations of motion, and encoding the forces of constraint in a bracket
of functions or a bilinear two-form. Once this is accomplished  the
constraints are satisfied automatically. The non-integrability of the constraint
distribution is then reflected in the failure of the bracket to satisfy the
Jacobi identity or in the failure of the bilinear two-form to be closed.

\subsubsection*{The Standard Almost Symplectic Formulation}
\label{S:The_Standard_Almost_Symplectic_Formulation}

In \cite{BS93} the Lagrange multipliers are eliminated from the equations of
motion 
(\ref{E:Eqns_of_motion_Ham_form_intrinsic})  using theorem \ref{T:Bates_Snyaticki}. Since $\mbox{$X_{\textup{nh}}^{\M}$}\in \C$ and the constraint forms $\tau^*\epsilon^i$ vanish along $\C$,
then the vector
field  $\mbox{$X_{\textup{nh}}^{\M}$}$ is uniquely determined by the equation
\begin{equation}
\label{E:Motion_Intrinsic_in_terms_of_Omega_C}
	{\bf i}_{\mbox{$X_{\textup{nh}}^{\M}$}} \, \Omega_\C = d\Ham_\C ,
\end{equation}
where $d\Ham_\C$ denotes the (point-wise) restriction of $d\Ham$ to $\C$. 
The latter equation really resembles the structure of a Hamiltonian system
except that 
$\Omega_\C$ is \emph{not} a two-form on $\M$.

\subsubsection*{The Standard Almost Poisson Formulation }
\label{S:AlmostPoissonFormulation}

We give a definition of the Almost Poisson Formulation for nonholonomic systems
that is convenient for our purposes. This formulation follows the general framework of Dirac brackets for systems with constraints and appears in the context of nonholonomic systems in \cite{Ibort99Dirac}. In   \cite{Cantrijn99} the authors show that the bracket so obtained coincides with that of \cite{vanderschaft1994, Marle1998}. It is shown in  \cite{vanderschaft1994} that this bracket is bilinear, anti-symmetric and satisfies Leibniz rule,  but it satisfies the Jacobi identity if and only if the constraints are holonomic. Thus the name \emph{almost} Poisson.

Let $\mathcal{P}:T_\M (T^*Q)\rightarrow \C$ be the projector associated to the decomposition
$T_\M (T^*Q)=\C \oplus \C^{\Omega_Q}$.

\begin{proposition}
\label{P:Important_for_Almost_Poisson_Bracket_Definition}
Let $f\in C^{\infty}(\M)$ and let $\bar f \in C^{\infty}(T^*Q)$ be an arbitrary smooth extension of $f$.
Let $X_{\bar f}$ be the free Hamiltonian vector field defined by ${\bf i}_{X_{\bar f}}\Omega_Q=d\bar f$.
Let $X_f^\C$ denote the unique vector field on $\M$ with values in $\C$ defined by the equation
\begin{equation}
\label{E:def_X_f^C}
{\bf i}_{X_f^\C}\; \Omega_\C=(d f)_\C,
\end{equation}
where $\Omega_\C$ and $(d f)_\C$ denote respectively the point-wise restriction of $\Omega_Q$ and
$df$ to $\C$. Then, along $\M$, we have $X_f^\C=\mathcal{P} X_{\bar f}$.
\end{proposition}
\begin{proof}
Let $m\in \M$. Applying $\mathcal{P}^*_m$ to both sides of ${\bf i}_{X_{\bar f}(m)}(\Omega_Q)_m=d\bar f(m)$ and pairing with an arbitrary $v_m\in T_m(T^*Q)$ gives, 
\begin{equation*}
\langle \mathcal{P}^*_m \; {\bf i}_{X_{\bar f}(m)}(\Omega_Q)_m \; , \; v_m \rangle =  \langle \mathcal{P}^*_m \; d\bar f(m) \; , \; v_m \rangle.
\end{equation*}
Since $\mathcal{P}_m$ is associated with a symplectic decomposition, we have
\begin{eqnarray*}
\langle \mathcal{P}^*_m \; {\bf i}_{X_{\bar f}(m)}(\Omega_Q)_m \; , \; v_m \rangle &=&  \langle {\bf i}_{X_{\bar f}(m)}(\Omega_Q)_m \; , \;  \mathcal{P}_m \; v_m \rangle \; = \;  (\Omega_Q)_m(X_{\bar f}(m)\; , \; \mathcal{P}_m \; v_m ) \\
&=& (\Omega_Q)_m(\mathcal{P}_m \; X_{\bar f}(m)\; , \;  v_m )\; = \; \langle {\bf i}_{\mathcal{P}_m \; X_{\bar f}(m)}(\Omega_Q)_m \; , \;  v_m \rangle.
\end{eqnarray*}
 It therefore follows that 
 \begin{equation}
 \label{E:auxinpropositionPoissonBracket}
 {\bf i}_{\mathcal{P}_m  X_{\bar f}(m)}(\Omega_Q)_m = \mathcal{P}^*_m \; d\bar f(m).
\end{equation}

By definition of $\mathcal{P}^*_m$, we have $\mathcal{P}^*_m \; d\bar f(m) = (d\bar f(m))_{\C_m}$.
 Since $\bar f$ is an extension of $f$,  the restriction of $d\bar f$ to $T\M$ agrees with $df$.
 In particular, the same is true about restriction to $\C_m \subset T_m\M$. Therefore $(d\bar f(m))_{\C_m}=(d f(m))_{\C_m}$. Moreover, since $\mathcal{P}_m  X_{\bar f}(m) \in \C_m$ we can replace $(\Omega_Q)_m$ by $(\Omega_{\C})_m$  in (\ref{E:auxinpropositionPoissonBracket}) to get
 \begin{equation*}
  {\bf i}_{\mathcal{P}_m  X_{\bar f}(m)}(\Omega_{\C})_m =(d f(m))_{\C_m}.
\end{equation*}
 By non-degeneracy of $(\Omega_{\C})_m$ we get $\mathcal{P}_m  X_{\bar f}(m)=X_f^\C(m)$ as required.
\end{proof}

In view of (\ref{E:Motion_Intrinsic_in_terms_of_Omega_C}) and the above proposition, it follows  that the 
nonholonomic vector field $\mbox{$X_{\textup{nh}}^{\M}$}$ satisfies
$\mbox{$X_{\textup{nh}}^{\M}$}=\mathcal{P}X_\Ham$, where $X_\Ham$ is the free Hamiltonian vector field defined by ${\bf i}_{X_\Ham}\Omega_Q=d\Ham$ and the equality makes sense on $\M$.

Let $f\in C^{\infty}(\M)$ and let $\bar f \in C^{\infty}(T^*Q)$ be an arbitrary smooth extension of $f$.
For $m\in \M$ we have:
\begin{eqnarray}
\nonumber
\mbox{$X_{\textup{nh}}^{\M}$}(f)(m)&=& \langle d\bar f (m) , \mathcal{P}_mX_\Ham(m) \rangle =(\Omega_Q)_m(X_{\bar f}(m),\mathcal{P}_mX_\Ham(m) )\\
&=& (\Omega_Q)_m(\mathcal{P}_mX_{\bar f}(m),\mathcal{P}_mX_\Ham(m) ), 
\label{E:Motivation_for_bracket_definition}
\end{eqnarray}
where the last identity follows from the fact that the projector $\mathcal{P}_m$ is associated to the symplectic 
decomposition $T_m (T^*Q)=\C_m \oplus \C_m^{\Omega_Q}$.
Inspired by this calculation we define the following bracket of functions  $f_1,f_2\in C^\infty(\M)$:
\begin{equation}
\label{E:bracket_definition}
\{f_1,f_2\}_\M(m):=(\Omega_Q)_m(\mathcal{P}_mX_{\bar f_1}(m),\mathcal{P}_mX_{\bar f_2}(m) )=\langle d f_1(m), \mathcal{P}_m X_{\bar f_2}(m) \rangle,
\end{equation}
where $\bar f_1, \bar f_2 \in C^{\infty}(T^*Q)$ are arbitrary smooth extensions of $f_1, f_2$.  The value of the bracket is independent of the extensions by proposition \ref{P:Important_for_Almost_Poisson_Bracket_Definition}. We will  refer to the above bracket as the  \emph{standard}  nonholonomic bracket  to distinguish it
from the \emph{affine} nonholonomic bracket to be introduced in section \ref{S:Affine_brackets}.

Associated to every function $f_1\in C^\infty(\M)$ we define its (almost) Hamiltonian vector field on $\M$,
denoted $X_{f_1}^\M$, by the rule
\begin{eqnarray*}
X_{f_1}^\M(f_2)(m)=\{f_2,f_1 \}_\M(m) \qquad \mbox{for all} \qquad f_2\in C^\infty(\M). 
\end{eqnarray*}
Denote by  $\Ham_\M$ is the restriction of $\Ham$ to $\M$. The equations of motion for the nonholonomic system can be written in terms of the standard nonholonomic bracket as:
\begin{eqnarray*}
\mbox{$X_{\textup{nh}}^{\M}$}(f)(m)=X_{\Ham_\M}^\M(f)= \{f,\Ham_\M \}_\M(m) \qquad \mbox{for all} \qquad f \in C^\infty(\M).
\end{eqnarray*}

\section{Reduction of Standard Nonholonomic Brackets}
\label{S:Reduction_of_Almost_Poisson_Brackets}

We present a basic scheme of reduction of nonholonomic systems, by reducing  the nonholonomic standard bracket. Our discussion is global and intrinsic and will be
extended for affine brackets ahead.
  A similar discussion is outlined in \cite{Marle1998}. In \cite{KoonMarsden2}  local  expressions are given for the reduced bracket and the link with  
 Lagrangian reduction is made.
 See also \cite{Naranjo2007} for other specialized cases of intrinsic reduction of almost Poisson brackets for nonholonomic systems on Lie groups.  
\begin{definition}
\label{D:Symmetry_group}
Let $H$ be a  Lie group that defines an action $\Phi:H\times Q\rightarrow Q$. We say that $H$ is a symmetry of our nonholonomic system if $\Phi$ lifts to a free and proper action on $TQ$ that 
leaves the constraint distribution $\D\subset TQ$ and the Lagrangian $\Lag:TQ\rightarrow \R$ invariant. 
\end{definition}

Suppose  that $H$ is a symmetry group of our nonholonomic system and denote by $\Psi:H\times T^*Q\rightarrow T^*Q$ the cotangent lift of $\Phi$. By the definition of cotangent lift, it follows that  $\Psi$ leaves the constraint submanifold, $\M$, and the Hamiltonian, $\Ham:T^*Q\rightarrow \R$,
invariant. The proof of the following proposition is left to the reader:

\begin{proposition}
The distribution $\F=\{ v\in T(T^*Q): \langle \tau^*\epsilon^i,  v\rangle =0\}\subset T(T^*Q)$ is invariant under
the lift of $\Psi$ to $T(T^*Q)$.
\end{proposition}

Since $\M$ is invariant under $\Psi$, the action naturally restricts to $\M$ and $T\M$ is invariant under the
tangent lift of $\Psi$. It follows from the above proposition that the restricted action to $\M$ preserves the
distribution $\C=T\M\cap \F$. Since $\Psi$ is the lift of a point transformation it is symplectic, so it
also preserves $\C^{\Omega_Q}$. As a consequence, it follows that for all $m\in \M$ and $h\in H$ the following diagram
commutes:
\begin{equation*}
\xymatrix{
{T_m(T^*Q)}
\ar @{->}[r]^-{T_m\Psi_h}
\ar @{->}[d]_-{\mathcal{P}_m}
&
{T_{\Psi_h(m)}(T^*Q)}
\ar[d]^-{\mathcal{P}_{\Psi_h(m)}}
\\
{T_m(T^*Q)}\ar @{->}[r]_-{T_m\Psi_h}
&
{T_{\Psi_h(m)}(T^*Q)}
}
\end{equation*}
It is now routine to check that $\Psi$ preserves the nonholonomic bracket. That is, for functions $f_1,f_2\in C^\infty(\M)$, and $h\in H$, we have
\begin{equation}
\label{E:invariant_bracket}
\{ f_1\circ \Psi_h,f_2\circ \Psi_h\}_\M=\{f_1,f_2\}_\M\circ \Psi_h.
\end{equation}
Let the smooth manifold $\mathcal{R}=\M/H$ denote the reduced space and $\pi: \M\rightarrow \mathcal{R}$ denote the orbit projection.
In view of (\ref{E:invariant_bracket}) the following \emph{reduced standard nonholonomic bracket} for functions  $F_1,F_2\in \mathcal{R}$ is well defined:
\begin{equation}
\label{E:def_of_reduced_bracket}
\{F_1,F_2\}_\mathcal{R}(\pi(m)):=\{F_1\circ \pi, F_2\circ \pi \}_\M(m).
\end{equation}
For a function $F_1\in C^\infty(\mathcal{R})$ we define its (almost) Hamiltonian vector field, $X_{F_1}^\mathcal{R}$ on $\mathcal{R}$ by the rule,
\begin{equation*}
X_{F_1}^\mathcal{R}(F_2)= \{F_2,F_1\}_\mathcal{R}, \qquad \mbox{for all} \qquad F_2\in C^\infty(\mathcal{R}).
\end{equation*}

Since the Hamiltonian $\Ham$ (and hence also the constrained Hamiltonian $\Ham_\M$) are
invariant, it follows that the nonholonomic vector field $\mbox{$X_{\textup{nh}}^{\M}$}=X_{\Ham_\M}^\M$ pushes forward by $\pi$ to 
the reduced nonholonomic vector field  $\mbox{$X_{\textup{nh}}^{\mathcal{R}}$}$ 
satisfying
\begin{equation*}
X_{\textup{nh}}^{\mathcal{R}}(F)=X_{\Ham_\mathcal{R}}^\mathcal{R}(F)=\{F,\Ham_\mathcal{R} \}_\mathcal{R} \qquad \mbox{for all} \qquad F \in C^\infty(\mathcal{R}),
\end{equation*}
where the reduced Hamiltonian $\Ham_\mathcal{R}\in C^\infty(\mathcal{R})$ is  uniquely defined
by the condition $\Ham_\M=\Ham_\mathcal{R}\circ \pi$. Summarizing, we have 
\begin{theorem}
\label{T:Almost_Poisson_Reduction} Suppose that the Lie group $H$ is a symmetry group for our nonholonomic system.
Then
\begin{enumerate}
\item The lifted action $\Psi$ on $T^*Q$ preserves the standard nonholonomic bracket $\{\cdot,\cdot \}_\M$ in the 
sense of (\ref{E:invariant_bracket}). 
\item The smooth reduced manifold $\mathcal{R}=\M/G$ is equipped
with a \emph{reduced standard nonholonomic  bracket}, $\{\cdot , \cdot \}_\mathcal{R}$, uniquely characterized by equation (\ref{E:def_of_reduced_bracket}).
\item The nonholonomic vector field, $\mbox{$X_{\textup{nh}}^{\M}$}$, is $\pi$-related to the (almost) Hamiltonian
vector field $X_{\Ham_\mathcal{R}}^\mathcal{R}$ associated to the reduced Hamiltonian $\Ham_\mathcal{R}$.
\end{enumerate}
\end{theorem}

It is easily shown that the reduced standard nonholonomic  bracket, $\{\cdot , \cdot \}_\mathcal{R}$, is almost Poisson.
A very interesting question is whether it in fact  satisfies the Jacobi identity thus yielding a direct Hamiltonization of the problem. It is shown in \cite{Naranjo2007} that this is indeed the case for the Suslov problem and the Chaplygin sleigh. A less stringent condition is the existence of a strictly positive function $\mu:Q/H \subset \mathcal{R}\rightarrow \R$ such that the new bracket of functions 
\begin{equation*}
\{F_1,F_2 \}_\mathcal{R}^\mu:=\mu \{F_1,F_2 \}_\mathcal{R},
\end{equation*}
satisfies the Jacobi identity. In this case we call $\mu$ a \emph{conformal factor} and the reduced equations
can be written in Hamiltonian form after the time rescaling $dt=\mu d\tau$.
Hamiltonization in this way is more likely to be accomplished if the   reduced space is low-dimensional and
 should not be expected in general. A necessary condition  is that the \emph{characteristic distribution}, $\mathcal{U}\subset T\mathcal{R}$, of the reduced standard bracket, defined by
 \begin{equation*}
\mathcal{U}_z=\{X^\mathcal{R}_F(z) \;:\; F\in C^\infty(\mathcal{R}) \}\subset T_z\mathcal{R}, 
\end{equation*}
  is integrable. This is a simple consequence of the symplectic stratification theorem for Poisson manifolds. In section  \ref{S:Chaplygin_Sphere} we will show that the characteristic distribution of the reduced standard bracket for the Chaplygin sphere is non-integrable, and thus, a conformal factor that renders this bracket Hamiltonian cannot exist. As mentioned
in the introduction, this example lead us to explore the more general concept of affine almost Poisson brackets and their reduction.

\section{Affine Almost Poisson Brackets}
\label{S:Affine_brackets}

In this section we will construct the affine almost Poisson brackets.
We begin by introducing the  notion of an \emph{Affine Almost Symplectic Structure} for 
a nonholonomic system. The idea is that the equations of motion for our nonholonomic system (\ref{E:Eqns_of_motion_Ham_form_intrinsic})
are equivalently written as
\begin{equation}
\label{E:Eqns_of_motion_Ham_form_intrinsic_Non_Standard}
{\bf i}_{\mbox{$X_{\textup{nh}}^{\M}$}} \,\iota^* (\Omega_Q + \Omega_0) = \iota^*(d\Ham +\lambda_i  \tau^* \epsilon^i),
\end{equation}
where $\Omega_0$ is any two-form on $T^*Q$ satisfying ${\bf i}_{\mbox{$X_{\textup{nh}}^{\M}$}} \,\iota^*  \Omega_0 =0$. 
We will require  the form $\tilde \Omega_Q=\Omega_Q+\Omega_0$ to be non-degenerate but
we will not ask for it to be closed. After all, the closeness of the canonical two-form $\Omega_Q$ was
never used in the construction of the standard nonholonomic bracket.

For our motivating example, the Chaplygin sphere, it is through the reduction of an  affine almost symplectic structure and a time rescaling that the system can be Hamiltonized.

\subsection*{The Affine Almost Symplectic Formulation}
\label{S:Affine_Almost_Symp_Formul}

We begin by giving our working definition of an affine almost symplectic structure.

\begin{definition}[Affine Almost Symplectic Structure] 
A nontrivial two-form $\Omega_0$ on $T^*Q$ defines an Affine Almost Symplectic Structure, $\tilde \Omega_Q:= \Omega_Q+\Omega_0$, for our nonholonomic
system if the following conditions hold:
\begin{enumerate}
\item ${\bf i}_{X_\Ham} \, \Omega_0= 0 $.
\item The form $\Omega_0$ is \emph{semi-basic} in the sense that it vanishes on vertical vectors. That is, if $v$ is a tangent vector to $T^*Q$ 
such that $\tau_*v=0$, with $\tau:T^*Q\rightarrow Q$ denoting the canonical projection, then ${\bf i}_{v} \, \Omega_0= 0$. \end{enumerate}
\end{definition}
 
 Condition 1 means that the form $\Omega_0$ ``does not see" the free Hamiltonian vector field $X_\Ham$. Condition  2 means that the form $\tilde \Omega_Q$ differs from the canonical form  by
 the addition of a magnetic type term. The following proposition shows that all the properties of $\Omega_Q$  that are relevant for the almost Hamiltonian formulation of nonholonomic systems are shared by 
$\tilde \Omega_Q$.   
 \begin{theorem}
\label{T:Properties_of_tilde_Omega} Let $\tilde \Omega_Q:= \Omega_Q+\Omega_0$ be an Affine Almost Symplectic Structure for our nonholonomic system.
The following statements are true:
\begin{enumerate}
\item The affine almost symplectic structure $\tilde \Omega_Q$ is non-degenerate. 
\item The point-wise restriction of $\tilde \Omega_Q$ to $\C$, denoted $\tilde \Omega_\C$, is non-degenerate.
\item The nonholonomic vector field $\mbox{$X_{\textup{nh}}^{\M}$}$ satisfies ${\bf i}_{\mbox{$X_{\textup{nh}}^{\M}$}}\iota^* \Omega_Q={\bf i}_{\mbox{$X_{\textup{nh}}^{\M}$}}\iota^* \tilde \Omega_Q$.
\item The symplectic complement of $\F$ with respect to $\Omega_Q$, denoted $\F^{\Omega_Q}$, equals
the symplectic complement of $\F$ with respect to $\tilde \Omega_Q$, denoted $\F^{\tilde \Omega_Q}$. That is $\F^{\Omega_Q}=\F^{\tilde \Omega_Q}$. 
\end{enumerate}
\end{theorem}
 \begin{proof}
 For $p_q\in T^*Q$, denote by $V_{p_q}\subset T_{p_q}(T^*Q)$ the subspace of vertical vectors, i.e. $V_{p_q}=\{v_{p_q}\in T_{p_q}(T^*Q): \tau_*v_{p_q}=0\}$. It is well known that $V_{p_q}$ is a Lagrangian
 subspace with respect to $\Omega_Q$, i.e. $V_{p_q}^{\Omega_Q}=V_{p_q}$.
 
 Let $w_{p_q}\in T_{p_q}(T^*Q)$ be such that $\tilde \Omega_Q(w_{p_q},v_{p_q})=0$ for all $v_{p_q}\in T_{p_q}(T^*Q)$. In particular, since $\Omega_0$ vanishes on vertical vectors, it follows
 that $\Omega_Q(w_{p_q},v_{p_q})=0$ for all $v_{p_q}\in V_{p_q}$, so $w_{p_q}\in V_{p_q}^{\Omega_Q}=V_{p_q}$. Using again that $\Omega_0$ vanishes on vertical vectors we get
 $\Omega_Q(w_{p_q},v_{p_q})=0$ for all $v_{p_q}\in T_{p_q}(T^*Q)$ which by non-degeneracy of $\Omega_Q$ implies $w_{p_q}=0$ and  we have proved 1.
 
Now,  for $m\in \M$ the intersection $V_m\cap \C_m$ is a Lagrangian subspace of $\C_m$ with respect to $(\Omega_\C)_m$. This follows from the  identity  $(V\cap \C)^{\Omega_Q}\cap \C=(V^{\Omega_Q}+ \C^{\Omega_Q})\cap \C=V\cap \C$.
 Repeating the argument in the above paragraph shows part 2.
 
 To prove 3 start by defining the vector field $\tilde X_\Ham$ by the equation $d\Ham= {\bf i}_{ \tilde X_{\Ham}}\tilde \Omega$. We claim that the vector fields  $\tilde X_\Ham$ and $X_\Ham$ are equal. Indeed, 
writing $\tilde X_\Ham = X_\Ham+Y_\Ham$, and since ${\bf i}_{ X_{\Ham}} \Omega_0=0$, we have\begin{equation*}
d\Ham ={\bf i}_{\tilde X_{\Ham}}\tilde \Omega_Q = {\bf i}_{ X_{\Ham}} \Omega_Q + {\bf i}_{Y_{\Ham}}\tilde \Omega_Q=d\Ham+ {\bf i}_{Y_{\Ham}}\tilde \Omega_Q.
\end{equation*}
It follows that ${\bf i}_{Y_{\Ham}}\tilde \Omega_Q=0$ which implies $Y_{\Ham}=0$ by non-degeneracy of $\tilde \Omega_Q$ shown above.

From the intrinsic form of the equations of motion 
(\ref{E:Eqns_of_motion_Ham_form_intrinsic}) we deduce that along $\M$ we can write
 $\mbox{$X_{\textup{nh}}^{\M}$}=X_\Ham+Z_\Ham=\tilde X_\Ham+Z_\Ham$, where
the \emph{constraint force vector field} $Z_\Ham\in \F^{\Omega_Q}$.  Since $V_{p_q}\subset \F_{p_q}$ for all $p_q\in T^*Q$, then $\F^{\Omega_Q}_{p_q} \subset V_{p_q}^{\Omega_Q}=V_{p_q}$ and  we conclude that $Z_\Ham$ is vertical. In view of the above observations and using properties (1) and (2) in the definition of the affine almost symplectic structure $\tilde \Omega_Q$ we find
\begin{equation*}
{\bf i}_{\mbox{$X_{\textup{nh}}^{\M}$}}\iota^* \Omega_Q =  {\bf i}_{\mbox{$X_{\textup{nh}}^{\M}$}}\iota^* \tilde \Omega_Q.
\end{equation*}

Finally, to prove 4, notice that since $\Omega_0$ vanishes on $V_{p_q}$, we have $V_{p_q}^{\tilde \Omega_Q}=V_{p_q}^{\Omega_Q}=V_{p_q}$.
 Now, by definition of $\F$ we have $V_{p_q}\subset \F_{p_q}$ for all $p_q\in T^*Q$. It follows that $\F^{\Omega_Q}_{p_q} \subset V_{p_q}^{\Omega_Q}=V_{p_q}$ and  $\F^{\tilde \Omega_Q}_{p_q} \subset V_{p_q}^{\tilde \Omega_Q}=V_{p_q}$. Since  the forms $\Omega_Q$ and $\tilde \Omega_Q$ agree when
 contracted with elements in $V_{p_q}$ the result follows.
 \end{proof}

 Point 3 in the above theorem shows that   starting from (\ref{E:Eqns_of_motion_Ham_form_intrinsic_Non_Standard}) and by a reasoning
 analogous to the discussion in section \ref{S:The_Standard_Almost_Symplectic_Formulation}, the vector field $\mbox{$X_{\textup{nh}}^{\M}$}$ is uniquely determined by the equation:
 \begin{equation}
 \label{E:Motion_Intrinsic_Affine_Structure_Rest_to_C}
	{\bf i}_{\mbox{$X_{\textup{nh}}^{\M}$}} \, \tilde \Omega_\C = d\Ham_\C .
\end{equation}

We finish the section with a small digression.  It is seen in the proof of point 3 in the above theorem, that the \emph{reaction force bundle} is the vector bundle over $T^*Q$ whose fibers are given by $\F^{\Omega_Q}$. Point 4 of the above  theorem shows that this bundle can also be written as the symplectic complement
of $\F$ with respect to $\tilde \Omega_Q$. This reinforces the idea that $\tilde \Omega_Q$ has indeed all of the relevant properties  for the description of nonholonomic systems that $\Omega_Q$ has.

\subsection*{The Affine Almost Poisson Formulation}
\label{SS:Affine_Almost_Poisson}

 Let $\tilde \Omega_Q=\Omega_Q +\Omega_0$ be an affine almost symplectic structure for our nonholonomic system.
In view of theorem \ref{T:Properties_of_tilde_Omega} we have the symplectic decomposition $T_m(T^*Q)=\C_m\oplus \C_m^{\tilde \Omega_Q}$. Let $\tilde {\mathcal{P}}:T_\M (T^*Q)\rightarrow \C$ be the projector associated to this  decomposition.
 Then, analogous to  proposition \ref{P:Important_for_Almost_Poisson_Bracket_Definition} one shows,
\begin{proposition}
\label{P:Important_for_Affine_Almost_Poisson_Bracket_Definition}
Let $f\in C^{\infty}(\M)$ and let $\bar f \in C^{\infty}(T^*Q)$ be an arbitrary smooth extension of $f$.
Let $\tilde X_{\bar f}$ be the  vector field on $T^*Q$ defined by ${\bf i}_{\tilde X_{\bar f}}\tilde \Omega_Q=d\bar f$.
Let $\tilde X_f^\C$ denote the vector field on $\M$ with values in $\C$ defined by the equation
\begin{equation*}
{\bf i}_{\tilde X_f^\C}\; \tilde \Omega_\C=(d f)_\C,
\end{equation*}
where $\tilde \Omega_\C$ and $(d f)_\C$ denote respectively the point-wise restriction of $\tilde \Omega_Q$ and
$df$ to $\C$. Then, along $\M$, we have $\tilde X_f^\C=\tilde{\mathcal{P}}\tilde X_{\bar f}$.
\end{proposition}
In view of (\ref{E:Motion_Intrinsic_Affine_Structure_Rest_to_C}) and as a consequence of the above
proposition, the 
nonholonomic vector field $\mbox{$X_{\textup{nh}}^{\M}$}$ satisfies $\mbox{$X_{\textup{nh}}^{\M}$}=\tilde{\mathcal{P}}\tilde X_{\Ham}$. 

 Analogous to (\ref{E:Motivation_for_bracket_definition}) we have:
\begin{eqnarray}
\nonumber
\mbox{$X_{\textup{nh}}^{\M}$}(f)(m)&=& \langle d\bar f (m) , \tilde{\mathcal{P}}_m\tilde X_\Ham(m) \rangle =(\tilde \Omega_Q)_m(\tilde X_{\bar f}(m),\tilde{\mathcal{P}}_m\tilde X_\Ham(m) ) \\
\label{E:Nonho_vector_field_in_terms_of_affine_bracket}
& =& (\tilde \Omega_Q)_m(\tilde{\mathcal{P}}_m\tilde X_{\bar f}(m),\tilde{\mathcal{P}}_m\tilde X_\Ham(m) ), 
\end{eqnarray}
where, for the last identity to hold,  we have used the fact that $\tilde{\mathcal{P}}$ is associated to a symplectic decomposition with respect to $\tilde \Omega_Q$. 
We now define the \emph{affine nonholonomic bracket}  for functions  $f_1,f_2\in C^\infty(\M)$:
\begin{equation*}
\{f_1,f_2 \tilde \}_\M(m):=(\tilde \Omega_Q)_m(\tilde{\mathcal{P}}_m\tilde X_{\bar f_1}(m),\tilde{\mathcal{P}}_m\tilde X_{\bar f_2}(m) )=\langle d f_1(m) ,   \tilde{\mathcal{P}}_m\tilde X_{\bar f_2}(m)\rangle, 
\end{equation*}
for arbitrary extensions $\bar f_1, \bar f_2 \in C^{\infty}(T^*Q)$ of $f_1,f_2$. The affine nonholonomic bracket is well defined in view of proposition \ref{P:Important_for_Affine_Almost_Poisson_Bracket_Definition}.

 Associated to every function $f_1\in C^\infty(\M)$ we define its (almost) affine Hamiltonian vector field on $\M$, denoted $\tilde X_{f_1}^\M$, by the rule
\begin{eqnarray*}
\tilde X_{f_1}^\M(f_2)(m)=\{f_2,f_1 \tilde{\}}_\M(m) \qquad \mbox{for all} \qquad f_2\in C^\infty(\M). 
\end{eqnarray*}

Notice that in general
$\{f_1,f_2\}_\M\neq \{f_1,f_2 \tilde \}_\M$ and consequently, for a general $f\in C^\infty(\M)$, the vector fields $X_{f}^\M$ and $\tilde X_f^\M$
are different. However, in view of (\ref{E:Nonho_vector_field_in_terms_of_affine_bracket}) we see that
\begin{eqnarray*}
\mbox{$X_{\textup{nh}}^{\M}$}(f)(m)=   X_{\Ham_\M}^\M(f)(m)=\tilde X_{\Ham_\M}^\M(f)(m)=\{f,\Ham_\M \tilde{\}}_\M(m), \qquad \mbox{for all} \qquad f\in C^\infty(\M), 
\end{eqnarray*}
so we can also
write the equations of motion for the nonholonomic system in terms of the affine bracket.

\section{Reduction of Affine Nonholonomic Brackets}
\label{S:Reduction_of_Affine_Nonholonomic_Brackets}

The reduction scheme presented in section \ref{S:Reduction_of_Almost_Poisson_Brackets} and summarized in theorem \ref{T:Almost_Poisson_Reduction} holds for affine brackets but we need to ask that the lifted action  to $T^*Q$ leaves the affine term $\Omega_0$ invariant.

\begin{theorem}
\label{T:Almost_Poisson_Reduction_Affine} Suppose that  $H$ is a symmetry group for our nonholonomic system.
Suppose, in addition, that there is an affine almost symplectic form $\tilde \Omega_Q=\Omega_Q + \Omega_0$ for our nonholonomic system and that the form $\Omega_0$ is invariant under the cotangent lifted action $\Psi:H\times T^*Q\rightarrow T^*Q$. Then
\begin{enumerate}
\item $\Psi$ preserves the affine bracket $\{\cdot,\cdot {\tilde \}}_\M$ in the 
sense that for all $f_1,f_2\in C^\infty(\M)$ and $h\in H$ we have
\begin{equation}
\{f_1\circ \Psi_h, f_2\circ \Psi_h \tilde{\}}_\M = \{f_1, f_2 \tilde{\}}_\M\circ \Psi_h.
\end{equation}
\item The smooth reduced manifold $\mathcal{R}=\M/G$ is equipped
with a \emph{reduced affine nonholonomic  bracket}, $\{\cdot , \cdot \tilde{\}}_\mathcal{R}$, uniquely determined by 
the relation
\begin{equation*}
\{F_1,F_2 \tilde{\}}_\mathcal{R}(\pi(m))=\{F_1\circ \pi ,F_2 \circ \pi \tilde{\}}_\mathcal{M}(m), \qquad \mbox{for}
\qquad F_1,F_2 \in C^\infty(\mathcal{R}), \;\; m\in \M,
\end{equation*}
and where $\pi:\M\rightarrow \mathcal{R}$ denotes the orbit projection.
\item
The nonholonomic vector field, $\mbox{$X_{\textup{nh}}^{\M}$}$, is $\pi$-related to the (almost) Hamiltonian
vector field $\tilde X_{\Ham_\mathcal{R}}^\mathcal{R}$ associated to the reduced Hamiltonian $\Ham_\mathcal{R}$, that is uniquely determined by the condition $\Ham_\M=\Ham_\mathcal{R}\circ \pi$.
\end{enumerate}  
\end{theorem}

In point 3 of the above theorem, and in what follows, we have denoted by $\tilde X_F^\mathcal{R}$ the (almost) Hamiltonian vector field corresponding to $F\in C^\infty(\mathcal{R})$ defined in terms of the reduced affine nonholonomic bracket $\{\cdot, \cdot \tilde{\}}_\mathcal{R}$.

Just as in the reduction of standard nonholonomic brackets discussed in section \ref{S:Reduction_of_Almost_Poisson_Brackets}, it is straightforward  to check that the  reduced affine nonholonomic  bracket, $\{\cdot , \cdot \tilde{\}}_\mathcal{R}$, is an almost Poisson bracket. The question
remains of whether it in fact satisfies the Jacobi identity, maybe after multiplication by a conformal factor.
We will show  that this is indeed the case for the Chaplygin sphere.

\section{Example: The Chaplygin Sphere}
\label{S:Chaplygin_Sphere}

We will illustrate the need for  affine almost Poisson brackets  by working out the almost Poisson reduction 
of the celebrated Chaplygin sphere problem with the standard nonholonomic bracket and with an affine bracket.
We will show how Hamiltonization after reduction can only be achieved by starting out with an affine bracket.

The Chaplygin sphere problem concerns the motion of an inhomogeneous sphere whose center 
of mass coincides with its geometric center that rolls without slipping on the plane. The configuration
of the system is $Q=SO(3)\times \R^2$. An element $q\in Q$ will be denoted by $q=(g;(x,y))$. Here $(x,y)$ give the cartesian coordinates of the center of the sphere on the plane, and 
$g\in SO(3)$ specifies the orientation of the sphere by relating, at any given time, a fixed space frame to a 
moving body frame that will be assumed to be aligned with the principal axes of inertia of the sphere.  
These two frames define respectively the so-called \emph{space} and \emph{body} coordinates.

The key aspect that renders the study of the Chaplygin sphere interesting is that the kinetic energy writes
naturally in terms of body coordinates while the constraints are naturally written in space coordinates. This is related to the fact that the problem is an LR system in the sense of \cite{Veselova} when considered on 
the \emph{direct} product Lie group $Q=SO(3)\times \R^2$. As a consequence, we will be forced to work with 
both the \emph{right} and the \emph{left} trivializations for $SO(3)$. Throughout the section we will 
constantly write formulas with respect to both coordinate systems.

\subsection*{A Moving Frame Approach}

To avoid working with Euler angles or other local coordinates for $SO(3)$ we will use moving frames to describe the system globally.

Identify the Lie algebra $\g=\so(3)$ with $\R^3$ using the \emph{hat-map},
\begin{equation*}
{\bf v}=(v_1,v_2,v_3)\mapsto \hat {\bf v}= \left ( \begin{array}{ccc} 0& -v_3 & v_2 \\ v_3 & 0 & -v_1 \\ -v_2 &v_1 & 0 \end{array} \right ).
\end{equation*}
The above is a Lie algebra isomorphism with the commutator in $\R^3$ being  the usual vector product.

Let $\{e_1, e_2, e_3 \}$ be the canonical basis for the Lie algebra $\g=\R^3$. The rolling takes place
on a plane parallel to $\{e_1,e_2\}$  and normal to $e_3$. The coordinates of the contact point  on the table with respect to the basis $\{e_1,e_2\}$ are  $(x,y)\in \R^2$.

The moving frame $\{X_1^\textup{right}(g),X_2^\textup{right}(g),X_3^\textup{right}(g)\}$ that forms a basis for $T_gSO(3)$ is obtained by right translation of $\{e_1, e_2, e_3 \}$. The dual co-frame will be denoted
by $\{ \rho_1(g),\rho_2(g),\rho_3(g)\}$. Similarly, left  translation of $\{e_1, e_2, e_3 \}$ by $g\in SO(3)$ defines the moving frame $\{X_1^\textup{left}(g),X_2^\textup{left}(g),X_3^\textup{left}(g)\}$ that forms a basis for $T_gSO(3)$ whose dual co-frame will be denoted
by $\{\lambda_1(g),\lambda_2(g),\lambda_3(g)\}$.

For a tangent vector $v_g\in T_gSO(3)$, the corresponding \emph{angular velocity in space coordinates} is the  element $\mathbf{ \omega}^s\in \g=\R^3$,  obtained by right trivialization. Its components are 
defined by
\begin{equation*}
v_g=\omega^s_1X_1^\textup{right}(g)+\omega^s_2X_2^\textup{right}(g)+\omega^s_3X_3^\textup{right}(g),
\end{equation*}
or, equivalently, by $\omega^s_i=\langle \rho_i(g) , v_g \rangle$. Analogously, the \emph{angular velocity in body coordinates} is the  element ${\bf \omega}^b\in \g=\R^3$,  obtained by left trivialization and whose components are 
given by $\omega^b_i=\langle \lambda_i(g) , v_g \rangle$.

These two vectors are related by the \emph{Adjoint map}: $\Ad_g:=T_e(L_g \circ R_{g^{-1}}):\g\rightarrow \g$.
We have $\mathbf{\omega}^b = \Ad_{g^{-1}}\mathbf{\omega}^s$. Define the coefficients $g_{ij}:SO(3)\rightarrow \R$ by
$\Ad_{g^{-1}}e_i=g_{ij}e_j$.
For ${\bf v} \in \R^3=\g$ we have $\Ad_g{\bf v}=g{\bf v}$. It follows that $\omega^b=g^{-1}\omega^s$ and that $g_{ij}$ equals the $i,j$ component ($i^{th}$ row, $j^{th}$ column) of the matrix $g\in SO(3)$. Since $g^{-1}=g^T$ we have $g_{ki}g_{kj}=\delta_{ij}$ and we can write
\begin{equation*}
\omega^b_j=g_{ij}\omega^s_i, \qquad \omega^s_i=g_{ij}\omega^b_j,
\end{equation*}
or equivalently, 
$\lambda_j(g)=g_{ij}\rho_i(g), \;\;\; \rho_i(g)=g_{ij}\lambda_j(g)$.

Denote by $c_{ikl}$  the structure constants of the Lie algebra defined by $[e_i,e_k]=c_{ikl}e_l$. We have  $c_{ikl}=0$ if two of the indices  are equal, $c_{ikl}=1$ 
if $(i, k, l)$ is a cyclic permutation of $(1, 2, 3)$ and 
$c_{ikl}=-1$ otherwise. The following proposition will be used in what follows.
\begin{proposition} 
\label{P:dif_R_i^j}
We have $dg_{ij}=c_{ikl}g_{lj}\rho_k=c_{lkj}g_{il}\lambda_k$. 
\end{proposition}
\begin{proof} Let $\{f_1,f_2,f_3\}$ be the dual basis to $\{e_1, e_2, e_3 \}$. We can write $g_{ij}=\langle f_j,\Ad_{g^{-1}}e_i \rangle$. Therefore,
\begin{eqnarray*}
dg_{ij}(X_k^\textup{right})&=&\frac{d}{dt} {\big |_{t=0}} \langle f_j,\Ad_{(\exp(e_kt)g)^{-1}}e_i \rangle =  \langle e^j,\Ad_{g^{-1}}[e_i,e_k] \rangle \\
&=& c_{ikl}\langle f_j,\Ad_{g^{-1}}e_l \rangle = c_{ikl}g_{lj}.
\end{eqnarray*}
For an arbitrary $v_g\in T_gSO(3)$, writing $v_g=\langle \rho_k, v_g\rangle X_k^\textup{right}$, and using the above equation shows the result. The proof of the other identity is analogous.
\end{proof}

\subsection*{The Chaplygin Sphere Problem}

At a given point $q=(g;(x,y))\in SO(3)\times \R^2$ we can use either $\{X_i^\textup{right}(g),
\partial_x,\partial_y\}$ or 
$\{X_i^\textup{left}(g),
\partial_x,\partial_y\}$  as basis for the tangent space $T_qQ$.
A tangent vector $v_q\in T_qQ$ can be written as 
\begin{eqnarray*}
v_q=\omega_i^sX_i^\textup{right}(g)+v_x\partial_x +v_y \partial_y  
=\omega_i^bX_i^\textup{left}(g)+v_x\partial_x +v_y \partial_y.
\end{eqnarray*}

The Lagrangian $\Lag:TQ\rightarrow \R$ is of pure kinetic energy and  in body coordinates is given by
\begin{eqnarray}
\label{E:Def_Lag_Body_Coords}
\Lag(v_q)=\frac{1}{2} ( \I {\bf \omega}^b) \cdot {\bf \omega}^b +\frac{1}{2}m(v_x^2 + v_y^2).
\end{eqnarray}
Here $\I$ is the inertia tensor which, under our assumption that the body frame is aligned with the principal axes of inertia of the body, is represented as a diagonal $3\times 3$ matrix whose positive entries, $I_i$, are the principal moments of inertia,   and ``\;$\cdot$\;" denotes the canonical scalar product on $\R^3$. 

The rolling constraints are:
\begin{equation}
\label{E:rolling_constraints}
v_x= r \omega^s_2,
\qquad
v_y = -r \omega^s_1,
\end{equation}
where $r$ is the radius of the sphere.
We have a three dimensional constraint distribution  $\D\subset TQ$ defined as the annihilator of  
\begin{eqnarray}
\label{E:Constraint_one_forms}
\epsilon_x= dx-r \rho_2, \qquad \epsilon_y= dy+r \rho_1.
\end{eqnarray}
A basis for $\D_q$ is $\{X_1^\textup{right}(g)-r\partial_y, X_2^\textup{right}(g)+r\partial_x, X_3^\textup{right}(g) \}$. The first two vector fields in this basis define rolling motions in the $x$ and $y$ directions with 
the accompanying rolling motion. The third one defines \emph{twisting} motion, where the sphere 
spins about the vertical $e_3$ axis but stays put in the table.

We now pass to the Hamiltonian formulation. The Lagrangian $\Lag$, being hyper-regular, defines a  Riemannian metric  on $Q$. For $q\in Q$, we have a natural isomorphism $\L_q:T_qQ\rightarrow T^*_qQ$, the Legendre transform of $\Lag$. The  Lagrangian then writes as $\Lag(v_q)=\frac{1}{2}\langle \L_q(v_q), v_q \rangle$ for $v_q\in T_qQ$,  where $\langle \cdot , \cdot \rangle$ is the duality pairing.

To compute explicitly the Legendre transform, notice that at the point $q=(g;(x,y))\in SO(3)\times \R^2$ we can use either $\{\rho_1(g),\rho_2(g),\rho_3(g),dx,dy\}$ or
$\{\lambda_1(g),\lambda_2(g),\lambda_3(g),dx,dy\}$ as basis for the cotangent space $T_q^*Q$.
A cotangent vector $\alpha_q\in T_q^*Q$, can then be written  as
\begin{eqnarray*}
\alpha_q=M^s_i\rho_i(g)+ p_xdx +p_ydy 
= M^b_i\lambda_i(g) + p_xdx +p_ydy.
\end{eqnarray*}
The vectors $\mathbf{M}^s$ and $\mathbf{M}^b$ with entries $M^s_i$ and $M^b_i$ are, respectively, the angular momentum of the ball with respect to its center expressed in space and body coordinates. They define elements in the dual Lie algebra  of $SO(3)$ and are related by $\mathbf{M}^b=\Ad_{g^{-1}}^*\mathbf{M}^s$. With the identification of $\R^3$ and $(\R^3)^*$ via the euclidean pairing, we can write
$\mathbf{M}^b=g^{-1}\mathbf{M}^s$.
 The vector $(p_x,p_y)$ is the linear momentum of the center of mass of the ball. We have thus defined 
 $(\mathbf{M}^s;p_x,p_y)$ and $(\mathbf{M}^b;p_x,p_y)$ as possible sets of coordinates for $T_q^*Q$. We will refer to them as space and body coordinates respectively.

The Legendre transformation maps the vector $v_q\in T_qQ$ to the covector $\alpha_q\in T^*_qQ$ according to the rule:
\begin{equation*}
\mathbf{M}^b=  \I {\bf \omega}^b, \qquad \mathbf{M}^s= g \I g^{-1}{\bf \omega}^s,\qquad p_x=mv_x, \qquad p_y=mv_y.
\end{equation*}

 The rolling constraints  write in the Hamiltonian formulation as,
$p_x=mr\omega^s_2, \;\; p_y=-mr\omega^s_1$,
with $\omega^s_1, \omega^s_2$ written in terms of $\mathbf{M}^b$ or $\mathbf{M}^s$ by means of the inverse Legendre transform. The constraint manifold $\M$ is defined as the set of $\alpha_q\in T^*Q$ such that these equations hold.

The Hamiltonian $\Ham:T^*Q\rightarrow \R$ is given in body coordinates by
\begin{eqnarray*}
\Ham(\alpha_q)= \frac{1}{2}\mathbf{M}^b\cdot (\I^{-1}\mathbf{M}^b) + \frac{1}{2m}(p_x^2+p_y^2).
\end{eqnarray*}

The canonical one-form on $T^*Q$ writes as\footnote{To simplify notation, and hoping that there is no danger of confusion,  we also denote by $\lambda_i, \; \rho_i$ the pull-backs $\tau^*\lambda_i, \; \tau^*\rho_i$ by the canonical projection $\tau:T^*Q\rightarrow Q$.} $\Theta_Q=M_i^b\lambda_i +p_xdx + p_ydy$ in the left trivialization, and as $\Theta_Q=M_i^s\rho_i +p_xdx + p_ydy$ in the right one.
Using Cartan's structure equations one has $d\lambda_1=-\lambda_2\wedge \lambda_3, \dots, \;\; d\rho_1=\rho_2\wedge \rho_3, \dots$ (cyclic), and we find the following expressions \ for the canonical two-form $\Omega_Q$,
\begin{eqnarray}
\nonumber
\Omega_Q=-d\Theta_Q
&=& \lambda_k\wedge dM^b_k+  M^b_1\lambda_2\wedge \lambda_3 + M_2^b \lambda_3\wedge \lambda_1 + M_3^b\lambda_1\wedge \lambda_2 + dx\wedge dp_x + dy \wedge dp_y \\
\label{E:Canonical_two_form}
&=& \rho_k \wedge dM^s_k-  M^s_1\rho_2\wedge \rho_3 - M_2^s \rho_3\wedge \rho_1 - M_3^s\rho_1\wedge \rho_2 + dx\wedge dp_x + dy \wedge dp_y,
\end{eqnarray}
 The free Hamiltonian vector field $X_\Ham$, defined by ${\bf i}_{X_\Ham}\Omega_Q=d\Ham$, is then shown to be given by\footnote{As with the previous remark, we ask the reader to interpret $X^\textup{left}_i, \partial_x, \partial_y$ in the right way. They  are tangent vectors to $T^*Q$ with no component in the $\partial_{M^b_i}, \partial_{p_x}, \partial_{p_y}$ directions, that push down to
the old $X^\textup{left}_i,\partial_x, \partial_y$ by $\tau:T^*Q\rightarrow Q$.}
\begin{eqnarray}
\label{E:Ham_vector_field_body}
X_\Ham=c_{kij}M_j^b\omega_k^b\partial_{M^b_i}+\omega^b_iX^\textup{left}_i+\frac{p_x}{m}\partial_x+\frac{p_y}{m}\partial_y,
\end{eqnarray}
in body coordinates. The above vector field defines the usual free rigid body equations for a body whose center of
mass lies on the $(x,y)$ plane. This is, of course,  trivial; in the absence of rolling constraint forces, the sphere is just a free rigid body. This type of motion will generally not  satisfy the constraints.

\subsection*{A Change of Coordinates on the Fibers}
\label{SS:Change_of_coordinates_on_fibers}

We  introduce a  change of coordinates on the fibers of $T^*Q$ inspired by 
\cite{vanderschaft1994}. See \cite{Koiller2002} for an interesting discussion of this kind of change of coordinates in the context of moving frames. In our particular case, the change of coordinates has a precise physical meaning; it corresponds to working  with the angular momentum with respect to the contact point,  instead of the angular momentum with respect to the center of mass.

Working in space coordinates, and since a basis for $\D_q$ is $\{X_1^\textup{right}-r\partial_y, X_2^\textup{right}+r\partial_x, X_3^\textup{right} \}$, following \cite{vanderschaft1994} we define at $T_q^*Q$ the new (space) fiber coordinates, $({\bf K}^s;m_x,m_y)$, in terms of our existing coordinates, $({\bf M}^s;p_x,p_y)$,   by
\begin{equation}
\label{E:Transformation_Rules}
K_1^s=M_1^s-rp_y, \qquad 
K_2^s=M_2^s+rp_x, \qquad 
K_3^s=M_3^s, \qquad
m_x=p_x, \qquad
m_y=p_y. 
\end{equation}
In body coordinates, we will consider the new fiber coordinates $({\bf K}^b;m_x,m_y)$ with 
${\bf K}^b:=g^{-1}{\bf K}^s$.

Along the constraint manifold $\M$ we
have $p_x=mr\omega^s_2$ and $p_y=-mr\omega^s_1$. Substituting these expressions into (\ref{E:Transformation_Rules}),  a short exercise in classical mechanics shows that  the vector ${\bf K}^s$ (resp. ${\bf K}^b$) is indeed  the
angular momentum of the ball with respect to the contact point written in space (resp. body) coordinates. 
Moreover, by writing $\mathbf{M}^s$ (resp.  $\mathbf{M}^b$) in terms of $\omega^s$ (resp. ${\bf \omega}^b$), one can easily derive the expressions,
\begin{eqnarray}
\label{E:Ks_in_terms_of_omegas}
{\bf K}^s&=&g \I g^{-1}{\bf \omega}^s +mr^2(\omega^s - (\omega^s\cdot e_3)e_3), \\
\label{E:Kb_in_terms_of_omegab}
{\bf K}^b&=& \I {\bf \omega}^b +mr^2(\omega^b - (\omega^b\cdot \gamma)\gamma),
\end{eqnarray}
where again ``$\;\cdot \;$" denotes the usual dot product in $\R^3$ and the Poisson vector, $\gamma\in \R^3$, is defined by 
\begin{equation*}
\gamma := g^{-1}e_3 \qquad \mbox{or,  in components,} \qquad  \gamma_i := g_{3i}.
\end{equation*}
 Physically, $\gamma$ is  the 
vertical unit vector expressed in body coordinates.  
For a Lie algebraic interpretation of $\gamma$ see \cite{Schneider}.

Expressions (\ref{E:Ks_in_terms_of_omegas}) and (\ref{E:Kb_in_terms_of_omegab}) can be inverted
to write $\omega^s$ (resp. $\omega^b$) in terms of ${\bf K}^s$ (resp. ${\bf K}^b$). To obtain  these  expressions explicitly
let $E$ denote the $3\times 3$ identity matrix and let $A:=(\I +mr^2E)$. Taking the dot product on both
sides of (\ref{E:Kb_in_terms_of_omegab}) with $A^{-1}\gamma$ yields, 
\begin{equation}
\label{E:OmegadotGamma}
\omega_3^s=\omega^b \cdot \gamma=\frac{\mathbf{K}^b\cdot A^{-1}\gamma}{Y(\gamma)},\qquad \mbox{with} \qquad Y(\gamma):=1-mr^2(\gamma\cdot A^{-1}\gamma).
\end{equation}
It follows that
\begin{equation}
\label{E:Omega_in_terms_of_K_and_Gamma}
\omega^b=A^{-1}\mathbf{K}^b+mr^2\left ( \frac{\mathbf{K}^b\cdot A^{-1}\gamma}{Y(\gamma)} \right ) A^{-1}\gamma,
\end{equation}
and
\begin{equation*}
\omega^s=gA^{-1}g^{-1}\mathbf{K}^s+mr^2\left ( \frac{(g^{-1}\mathbf{K}^s)\cdot A^{-1}\gamma}{Y(\gamma)} \right ) gA^{-1}\gamma.
\end{equation*}

The dimensionless quantity $Y(\gamma)$ will turn out to be very important in the context of Hamiltonization. 
It is seen to be strictly positive since $||\gamma||=1$ and all of the principal moments of inertia 
$I_i>0$.

Along the constraint submanifold $\M$, the above equations allow us to write, via the constraint equations $m_x=p_x=mr\omega^s_2, \;\; m_y=p_y=-mr\omega^s_1$, the quantities $m_x$ and $m_y$ as functions of $(g,{\bf K}^s)$ or $(g,{\bf K}^b)$. We can therefore use $((g;(x,y));,{\bf K}^s)$ or $((g;(x,y)),{\bf K}^b)$ as induced
coordinates for $\M$. The constrained Hamiltonian $\Ham_\M:= \Ham |_\M$ can be written in these
coordinates as
\begin{equation}
\label{E:constrained_Hamiltonian}
\Ham_\M=\frac{1}{2}\mathbf{K}^s \cdot \omega^s= \frac{1}{2}\mathbf{K}^b \cdot \omega^b.
\end{equation}

\subsection*{The $SE(2)$ Symmetry}
\label{E:the_SE(2)_symmetry}

Since the time evolution of the system is independent of the position of the origin and the orientation of the 
axes on the plane where the rolling takes place, it is natural to expect the left multiplication by elements in $H=SE(2)$  to be a symmetry group for the problem.

We represent $SE(2)$ as the subgroup of $GL_4(\R)$ consisting  of  matrices of the form
\begin{equation*}
\left ( \begin{array}{cc} \begin{array}{c} h \end{array} & \begin{array}{c} a \\ b \\ 0  \end{array} \\ \begin{array}{ccc}0& 0 & 0 \end{array} & 1 \end{array}  \right ), \qquad \mbox{where $h\in SO(3)$ is of the form} \qquad h=\left ( \begin{array}{cc} \tilde h & \begin{array}{c} 0 \\ 0 \end{array} \\ \begin{array}{cc} 0 & 0 \end{array} & 1 \end{array} \right ),
\end{equation*}
 with $\tilde h\in SO(2)$ and $a,b\in \R$. The generic element in $\se(2)$ will be denoted by $(h;a,b)$.
The action on an element $q=(g;x,y)\in Q$ is defined as:
\begin{equation*}
(h;a,b)\;:\; (g;x,y)\longrightarrow (hg; (x,y)\tilde h^t+(a,b))\in Q.
\end{equation*}

\begin{proposition} With the above definition of its  action on $Q$, $SE(2)$  is a symmetry group for
the Chaplygin sphere system.
\end{proposition}
\begin{proof}
We need to show that the lifted action to $TQ$ leaves both the constraints and the constraint distribution
invariant. Working with body coordinates, the lifted 
action to $TQ$ maps the tangent vector $v_q=(\omega^b; v_x,v_y)\in T_qQ$ according to
the rule
\begin{equation*}
\label{E:Tangent_lifted_action_body}
(h;a,b)\;:\; (\omega^b; v_x,v_y)\longrightarrow (\omega^b; (v_x,v_y)\tilde h^t)\in T_{(h;a,b)\cdot q}Q,
\end{equation*}
and it is immediate to check that the Lagrangian (\ref{E:Def_Lag_Body_Coords}) is invariant.
Working in space coordinates, the lifted action to $TQ$ maps the tangent vector $v_q=(\omega^s; v_x,v_y)\in T_qQ$ according to
the rule
\begin{equation}
\label{E:Tangent_lifted_action_space}
(h;a,b)\;:\; (\omega^s; v_x,v_y)\longrightarrow (\Ad_h\omega^s; (v_x,v_y)\tilde h^t)=(h\omega^s; (v_x,v_y)\tilde h^t)\in T_{(h;a,b)\cdot q}Q,
\end{equation}
and it is readily shown that the rolling constraints (\ref{E:rolling_constraints}) are invariant.
\end{proof}

According to theorem \ref{T:Almost_Poisson_Reduction}, the reduced space $\mathcal{R}:=\M/SE(2)$ is equipped with the reduced standard nonholonomic bracket $\{\cdot, \cdot \}_\mathcal{R}$, and the reduced
equations can be written with respect to this bracket. The properties of this bracket will be studied  further ahead.

We claim that associated with this symmetry there is a  conservation law: the vertical component of angular momentum about the center of mass of the sphere, $M_3^s$, (which agrees with vertical component of the sphere's angular momentum about the contact point, $K_3^s$) is constant throughout the motion. Contrary to usual holonomic mechanics, the relationship between symmetries and conservation laws is \emph{not} straightforward for nonholonomic systems, (see \cite{BKMM} for a thorough discussion). To show our claim we could follow \cite{BKMM} and compute the \emph{momentum equation} but we choose  to invoke instead the following 
nonholonomic version of 
Noether's theorem whose proof can be found in  \cite{ArnoldIII}:

\begin{theorem}\label{T:NoetherNonho} 
Let $H$ be a Lie group that acts on $Q$ with Lie algebra $\h$ and dual Lie algebra $\h^*$. For $\xi\in \h$, denote by $\xi_Q(q)\in T_qQ$ the infinitesimal generator of the action on $Q$, and by $\xi_{T^*Q}(\alpha_q)\in T_{\alpha_q}(T^*Q)$ the infinitesimal generator of the lifted action to $T^*Q$. If the Lagrangian $\Lag$ is invariant under the lifted action to $TQ$ and if $\xi_Q(q)\in \D_q$ for all $q\in Q$, then  the components of the momentum map, 
$\J_H:T^*Q\rightarrow \h^*$, defined by
\begin{equation*}
d\langle \J_H(\alpha_q), \xi \rangle = {\bf i}_{\xi_{T^*Q}(\alpha_q)}\Omega_Q,
\end{equation*}
are constant during the motion.
\end{theorem}
Notice that, contrary to the  definition of a symmetry group for a nonholonomic system given in section
\ref{S:Reduction_of_Almost_Poisson_Brackets},
in the above theorem we do \emph{not} require the constraint distribution to be invariant under the action. Notice as well that the momentum map, $\J_H:T^*Q\rightarrow \h^*$, always exists since we are working with a lifted action,  see \cite{MarsdenRatiubook}.

To apply the theorem we consider the \emph{twisting} action of $SO(2)$ on $Q$.  
Represent $SO(2)$ as the subgroup of $SO(3)$ consisting of matrices of the form,
\begin{equation*}
h=\left ( \begin{array}{cc} \tilde h & \begin{array}{c} 0 \\ 0 \end{array} \\ \begin{array}{cc} 0 & 0 \end{array} & 1 \end{array} \right ), \qquad \mbox{with} \qquad \tilde h\in SO(2).
\end{equation*}
The twisting $SO(2)$ action on $Q$ is defined by:
\begin{equation*}
h \;:\; (g;(x,y))\longrightarrow (hg; (x,y)), \qquad \mbox{for} \qquad (g;(x,y))\in Q.
\end{equation*}
The proof of the following proposition is left to the reader.
\begin{proposition}
The twisting $SO(2)$ action on $Q$ satisfies the hypothesis of theorem \ref{T:NoetherNonho} and
the associated conserved quantity is $M_3^s=K_3^s=\mathbf{M}^b\cdot \gamma= \mathbf{K}^b\cdot \gamma$.
\end{proposition}

After this digression on conserved quantities we come back to the discussion of the $SE(2)$ symmetry with the aim of performing the corresponding reduction as described in section \ref{S:Reduction_of_Almost_Poisson_Brackets}. In view of (\ref{E:Tangent_lifted_action_space}),
 the cotangent lift of the action to $T^*Q$ expressed in space coordinates maps the cotangent vector $ (\textbf{M}^s; p_x,p_y) \in T_q^*Q$
according to the rule
\begin{equation}
\label{E:cotangent_lifted_action_spaceM_variables}
(h;a,b)\;:\; (\textbf{M}^s; p_x,p_y)\longrightarrow (\Ad^*_{h^{-1}}\textbf{M}^s; (p_x,p_y)\tilde h^t)=(h\textbf{M}^s; (p_x,p_y)\tilde h^t)\in T^*_{(h;a,b)\cdot q}Q.
\end{equation}

 In view of (\ref{E:cotangent_lifted_action_spaceM_variables}) and the transformation rules (\ref{E:Transformation_Rules}), the cotangent lift of $SE(2)$ to $T^*Q$ is expressed in the new space fiber  coordinates
as the map acting on  $({\bf K}^s;m_x,m_y)\in T_q^*Q$ by
\begin{equation*}
(h;a,b)\;:\; ({\bf K}^s; m_x,m_y)\longrightarrow (\Ad^*_{h^{-1}}{\bf K}^s; (m_x,m_y)\tilde h^t)=(h{\bf K}^s; (m_x,m_y)\tilde h^t)\in T^*_{(h;a,b)\cdot q}Q.
\end{equation*}
Consequently, the action  in the new body coordinates  $({\bf K}^b;m_x,m_y)\in T_q^*Q$ is given by
\begin{equation*}
(h;a,b)\;:\; ({\bf K}^b;m_x,m_y)\longrightarrow ( {\bf K}^b; (m_x,m_y)\tilde h^t)\in T^*_{(h;a,b)\cdot q}Q.
\end{equation*}

Since the action preserves the constraint manifold $\M$, the restricted action on $\M$ 
is represented in the induced  body coordinates $((g;(x,y));{\bf K}^b)$ for $\M$ as
\begin{equation*}
(h;a,b)\;:\;((g;(x,y));{\bf K}^b) \longrightarrow ((hg; (x,y)\tilde h^t + (a,b)); {\bf K}^b).
\end{equation*}

Notice that the action leaves ${\bf K}^b$ and the components of the Poisson vector $\gamma$ (the third row of $g$) invariant. The latter are not independent as $||\gamma||=1$, but we can  use
$({\bf K}^b,\gamma)$ as redundant coordinates for the reduced space $\mathcal{R}:=\M/SE(2)\cong \R^3\times S^2$. In what follows we will represent $\mathcal{R}$ as the embedded submanifold
in $\R^6=\{({\bf K}^b,\gamma)\;:\; {\bf K}^b\in \R^3, \; \gamma \in \R^3 \}$ defined by the condition $||\gamma||=1$.

 The reduced Hamiltonian $\Ham_\mathcal{R}:\mathcal{R}\rightarrow \R$ is given by
 \begin{equation}
 \label{E:Reduced_Hamiltonian}
\Ham_\mathcal{R}(\mathbf{K}^b,\gamma)=\frac{1}{2}\mathbf{K}^b\cdot \omega^b,
\end{equation}
 with $\omega^b$ given by (\ref{E:Omega_in_terms_of_K_and_Gamma}). Finally notice that the 
 constant of motion, $M_3^s=K_3^s={\bf K}^b\cdot \gamma$, is invariant under the action since it can be written in terms of the coordinates $({\bf K}^b,\gamma)$.

\subsection*{An Affine Almost Symplectic Structure for the Chaplygin Sphere}
\label{SS:Affine_Structure}

Let $\nu$ denote the dimensionless, bi-invariant volume form on $SO(3)$, oriented and scaled such that for the canonical vectors $e_1,e_2,e_3\in \R^3$ we have
$\nu(e_1,e_2,e_3)=1$.
Since $Q=SO(3)\times \R^2$,  $\nu$ naturally defines a three-form on   $Q$ that, via the  cotangent bundle projection, $\tau:T^*Q\rightarrow Q$,  pulls back to a basic three-form $\bar \nu$ on $T^*Q$.

Denote by $X_\Ham$ the free Hamiltonian vector field of the sphere  (\ref{E:Ham_vector_field_body}) and define the two-form $\Omega_0$ on $T^*Q$ by
\begin{equation}
\label{E:def_Omega_0_for_sphere}
\Omega_0:=-mr^2 \; {\bf i}_{X_\Ham}\bar \nu.
\end{equation}

\begin{proposition}
The two-form $\tilde \Omega_Q:= \Omega_Q + \Omega_0$ on $T^*Q$ defines an affine almost symplectic structure for the Chaplygin sphere problem.
\end{proposition}
\begin{proof}
It is clear that $\Omega_0$ is a semi-basic two-form since it was constructed by pulling back a form on $Q$ and then contracting with $X_\Ham$. It is also clear that
${\bf i}_{X_\Ham} \Omega_0=0$ so the two conditions in the definition of an affine almost symplectic
structure are satisfied.
\end{proof}

Therefore, according to the theory developed in section \ref{S:Affine_brackets}, there is an affine nonholonomic bracket, $\{\cdot, \cdot \tilde{\}}_\M$, associated with the affine almost symplectic structure $\tilde \Omega_Q$. The following proposition shows that the hypothesis in theorem \ref{T:Almost_Poisson_Reduction_Affine} are satisfied so we can  reduce the system in terms of an affine reduced nonholonomic bracket $\{\cdot, \cdot \tilde{\}}_\mathcal{R}$ on $\mathcal{R}=\M/SE(2)$.

\begin{proposition}
The form $\Omega_0$ defined above is invariant under the cotangent lift  of the $SE(2)$ action to $T^*Q$. 
\end{proposition}
\begin{proof}
For $\xi\in \se(2)$, denote by $\xi_{T^*Q}$ the infinitesimal generator of the action on $T^*Q$. We have
\begin{equation*}
\pounds_{\xi_{T^*Q}}\Omega_0=-mr^2\left ( {\bf i}_{[\xi_{T^*Q}, X_\Ham]}\bar \nu +  {\bf i}_{ X_\Ham}\pounds_{\xi_{T^*Q}}\bar \nu \right ).
\end{equation*}
Since the Hamiltonian $\Ham:T^*Q\rightarrow \R$ is invariant under the action we have $[\xi_{T^*Q}, X_\Ham]=0$. We also have $\pounds_{\xi_{T^*Q}}\bar \nu=0$ by invariance of $\bar \nu$. Thus, $\pounds_{\xi_{T^*Q}}\Omega_0=0$, and the result follows since $Q$ is connected.
\end{proof}

Since the three-form $\bar \nu$ is given by $\bar \nu=\lambda_1\wedge \lambda_2\wedge \lambda_3$ in body coordinates.   In view of (\ref{E:Ham_vector_field_body}) we find the following explicit formula for $\Omega_0$., 
\begin{eqnarray}
\label{E:Explicit_Omega_0}
\Omega_0=-mr^2(\omega^b_1\lambda_2\wedge \lambda_3 + \omega^b_2\lambda_3\wedge \lambda_1+\omega^b_3\lambda_1\wedge \lambda_2).
\end{eqnarray}

\subsection*{The Reduced Brackets}
\label{SS:Reduced_brackets}

We now study the reduced brackets $\{\cdot,\cdot \}_\mathcal{R}$, and $\{\cdot,\cdot \tilde{\}}_\mathcal{R}$, on the reduced 
5-dimensional space $\mathcal{R}\cong \R^3\times S^2$. We will show that the characteristic distribution corresponding to the reduced standard nonholonomic bracket, $\{\cdot, \cdot \}_\mathcal{R}$
 is non-integrable, being thus very different
from that of a true Poisson bracket. In contrast, we will  show that the reduced affine nonholonomic bracket, $\{\cdot, \cdot \tilde{\}}_\mathcal{R}$, satisfies the Jacobi identity after
multiplication by a conformal factor. This allows us to effectively write the reduced equations of motion 
in Hamiltonian form after a rescaling of time.

Recall that the reduced space $\mathcal{R}$ is represented as the embedded submanifold in 
$\R^6=\{({\bf K}^b,\gamma)\;:\; {\bf K}^b\in \R^3, \; \gamma \in \R^3 \}$ defined by the condition $||\gamma||=1$. We will give explicit formulae for the brackets in these coordinates but we first need to derive some  formulas.

We begin by writing the canonical symplectic form $\Omega_Q$ in terms of our new fiber coordinates.
In view of (\ref{E:Canonical_two_form}) and (\ref{E:Transformation_Rules}) we find,
\begin{eqnarray*}
\Omega_Q&=&\rho_i\wedge dK^s_i-  K^s_1\rho_2\wedge \rho_3 - K_2^s \rho_3\wedge \rho_1 - K_3^s\rho_1\wedge \rho_2 -rm_y\rho_2\wedge \rho_3+rm_x\rho_3\wedge \rho_1 \\ &&+ \;\; (\tau^* \epsilon_x)\wedge dm_x + (\tau^*\epsilon_y) \wedge dm_y,
\end{eqnarray*}
where $\epsilon_x, \; \epsilon_y$ are the constraint one-forms on $Q$ defined in
(\ref{E:Constraint_one_forms}). The pull-back of $\Omega_Q$ to $\M$ via the inclusion map $\iota:\M \hookrightarrow T^*Q$ can be expressed in the induced coordinates $((g;(x,y),\mathbf{K}^s)$ as,
\begin{eqnarray*}
\iota ^*\Omega_Q&=&\rho_i\wedge dK^s_i-  K^s_1\rho_2\wedge \rho_3 - K_2^s \rho_3\wedge \rho_1 - K_3^s\rho_1\wedge \rho_2 +mr^2(\omega^s_1\rho_2\wedge \rho_3+ \omega^s_2 \rho_3\wedge \rho_1) \\ && \qquad + mr( (\tau^*\epsilon_x)\wedge d\omega^s_2 - (\tau^*\epsilon_y) \wedge d\omega^s_1).
\end{eqnarray*}
Therefore, the restriction, $\Omega_\C$, of $\iota ^*\Omega_Q$ to the space $\C=T\M\cap \mbox{ann}\{\tau^*\epsilon_x, \tau^*\epsilon_y\} \subset T(T^*Q)$ is given by
\begin{equation}
 \label{E:Omega_C_in_space_coordinates}
\Omega_\C=\rho_i\wedge dK^s_i-  K^s_1\rho_2\wedge \rho_3 - K_2^s \rho_3\wedge \rho_1 - K_3^s\rho_1\wedge \rho_2 +mr^2(\omega^s_1\rho_2\wedge \rho_3+ \omega^s_2 \rho_3\wedge \rho_1).
\end{equation}
\begin{proposition}
The restricted forms  $\Omega_\C$ and $\tilde \Omega_\mathcal{C}$ are expressed in body coordinates as 
\begin{eqnarray}
\nonumber
\Omega_\C&=&\lambda_i\wedge dK^b_i+  K^b_1\lambda_2\wedge \lambda_3 + K_2^b \lambda_3\wedge \lambda_1 + K_3^s\lambda_1\wedge \lambda_2 +mr^2(\omega^b_3\lambda_1\wedge \lambda_2+ \omega^b_1 \lambda_2\wedge \lambda_3 + \omega^b_2\lambda_3 \wedge \lambda_1) \\
 \label{E:Omega_C_in_body_coordinates}
 && \qquad  -mr^2\omega_3^s(\gamma_3 \lambda_1\wedge \lambda_2 + \gamma_2\lambda_3\wedge \lambda_1 + \gamma_1\lambda_2\wedge \lambda_3), \\
 \nonumber
 \tilde \Omega_\C&=&\lambda_i\wedge dK^b_i+  K^b_1\lambda_2\wedge \lambda_3 + K_2^b \lambda_3\wedge \lambda_1 + K_3^s\lambda_1\wedge \lambda_2 \\
 \label{E:Omega_Tilde_C_in_body_coordinates}
 && \qquad  -mr^2\omega^s_3(\gamma_3 \lambda_1\wedge \lambda_2 + \gamma_2\lambda_3\wedge \lambda_1 + \gamma_1\lambda_2\wedge \lambda_3).
\end{eqnarray}
\end{proposition}
\begin{proof}
Adding and subtracting $mr^2(\omega^s_3\rho_1\wedge \rho_2)$ to (\ref{E:Omega_C_in_space_coordinates}) and using the identities $K_i^s=g_{ij}K_j^b$, $\omega_i^s=g_{ij}\omega_j^b$, and $\rho_i=g_{il}\lambda_l$, together with proposition \ref{P:dif_R_i^j}, we get,
\begin{eqnarray*}
\nonumber
\Omega_\C&=&g_{il}g_{ij}\lambda_l \wedge dK^b_j+ c_{rkj}g_{il}g_{ir}  K^b_j\lambda_l\wedge \lambda_k \\  && \qquad +
(g_{1j}g_{2l}g_{3r}+g_{2j}g_{3l}g_{1r}+g_{3j}g_{1l}g_{2r}) (-K_j^b+mr^2\omega_j^b) \lambda_l \wedge \lambda_r -mr^2 \omega_3^s ( g_{1j}g_{2l})\lambda_j\wedge \lambda_l .
\end{eqnarray*}
Since $g\in SO(3)$ we have, $g_{il}g_{ij}=\delta_{lj}, \; \det(g)=1$, and $g_{1j}g_{2l}-g_{2j}g_{1l}=c_{jlk}g_{3k}=c_{jlk}\gamma_{k}$. Using these identities in the above expression for $\Omega_\mathcal{C}$ gives (\ref{E:Omega_C_in_body_coordinates}). The proof of (\ref{E:Omega_Tilde_C_in_body_coordinates}) is immediate in view of (\ref{E:Explicit_Omega_0}).
\end{proof}

 Use $((g;x,y);{\bf K}^b)$ as coordinates for $\M$. Denote by $\hat \partial_{K^b_i}\in T\M$ the tangent vector obtained as a  derivation with respect to $K^s_i$  when considered
as a coordinate on $\M$. Similarly, denote by $\hat X_i^\textup{left}, \hat \partial_x, \hat \partial_y$ the tangent vectors to $\M$, with zero component in the directions of $\hat \partial_{K^b_j}$ that push forward to
the tangent vectors $ X_i^\textup{left}, \partial_x, \partial_y \in TQ$ by the composition $\tau \circ  \iota :\M \hookrightarrow T^*Q \rightarrow  Q$. The latter always exist since $\M$ is a vector bundle over $Q$. We claim that 
\begin{equation}
\label{E:Basis_for_C}
\mathcal{B}:=\{\hat X_i^\textup{left}+r(g_{2i}\hat \partial_x-g_{1i}\hat \partial_y), \;  \hat \partial_{K^b_i} \;  : \; i=1,2,3 \}
\end{equation}
is a basis for the subspace $\mathcal{C}=T\M\cap \mbox{ann}\{\tau^*\epsilon_x, \tau^*\epsilon_y\}$. 
The vectors are tangent to $\M$ by definition. Moreover,   they are annihilated by
$\tau^*\epsilon_x$ and  $\tau^*\epsilon_y$ since in body coordinates we
have 
\begin{equation*}
\epsilon_x=dx-r g_{2j}\lambda_j, \qquad \epsilon_y=dy+r g_{1j}\lambda_j.
\end{equation*}
It is also immediate to check that they are linearly independent and span $\C$.

\subsubsection*{The Reduced Standard Nonholonomic Bracket}

 The following proposition gives explicit formulae for the reduced standard bracket.

\begin{proposition} 
\label{P:Formulae_for_reduced_standard_bracket}
We have
\begin{equation*}
\{K^b_i,K^b_j\}_\mathcal{R}=-c_{ijl} \left (  K^b_l +mr^2(  \omega^b_l - \omega^s_3\gamma_l )\right ), \qquad \{K^b_i,\gamma_j \}_\mathcal{R}= -c_{ijl} \gamma_l,  \qquad \{\gamma_i,\gamma_j \}_\mathcal{R}=0.
\end{equation*}
\end{proposition}
\begin{remark} Notice that the quantities $\omega^s_3$ and $\omega^b$, appearing in the above formulae, can be expressed in terms of our coordinates   $({\bf K}^b,\gamma)$ via (\ref{E:OmegadotGamma}) and
(\ref{E:Omega_in_terms_of_K_and_Gamma}). Also notice, by a short calculation, that according to the above formulae, $||\gamma||^2$ is a Casimir function,  so they indeed define a bracket on $\mathcal{R}$.
\end{remark}
\begin{proof}
Use $((g;x,y);{\bf K}^b)$ as coordinates for $\M$. With the definition given for equation (\ref{E:def_X_f^C}),
we claim that 
 \begin{equation*}
X_{K^b_i}^\C=\hat X_i^{\textup{left}}+r(g_{2i} \hat \partial_x-g_{1i}\hat \partial_y)+c_{ijl}(K_l^b +mr^2( \omega^b_l - \omega^s_3\gamma_l ))\hat \partial_{K_j^b}, \qquad X_{\gamma_i}^\C=c_{ijl}\gamma_l \hat \partial_{K_j^b}.
\end{equation*}
To show the claim, first note that the above vector fields indeed lie on $\C$ as they are expressed in terms of the basis $\mathcal{B}$ defined in (\ref{E:Basis_for_C}).
Next, putting $\gamma_i=g_{3i}$ and using proposition \ref{P:dif_R_i^j} we find $d\gamma_i=-c_{ijl}\gamma_l\lambda_j$, and in view of (\ref{E:Omega_C_in_body_coordinates}) one verifies by a direct calculation that the above vector fields satisfy:
\begin{equation*}
{\bf i}_{X_{K^b_i}^\C}\Omega_\C=(dK^b_i)_\C, \qquad {\bf i}_{X_{\gamma_i}^\C}\Omega_\C=(d\gamma_i)_\C.
\end{equation*}
In addition,  by proposition \ref{P:Important_for_Almost_Poisson_Bracket_Definition}, the above vector fields equal $\mathcal{P} X_{K^b_i}$ and $\mathcal{P} X_{\gamma_i}$ respectively. Therefore, by  definition of  the standard nonholonomic bracket $\{\cdot ,\cdot \}_\M$  given in (\ref{E:bracket_definition}) we have, 
\begin{eqnarray*}
\{K^b_i,K^b_j\}_\mathcal{M}=-\langle dK^b_j,\mathcal{P} X_{K^b_i} \rangle,   \qquad    \{K^b_i,\gamma_j \}_\mathcal{M}=-\langle d\gamma_j, \mathcal{P} X_{K^b_i} \rangle,  \qquad \{\gamma_i,\gamma_j \}_\mathcal{M}=-\langle d\gamma_j,\mathcal{P} X_{\gamma_i} \rangle.  
\end{eqnarray*}
The result now follows by computing the above pairings explicitly, and then using invariance of $\mathbf{K}^b$ and $\gamma$ and the definition of the reduced standard
nonholonomic bracket $\{\cdot ,\cdot \}_\mathcal{R}$.
\end{proof}

In view of the above proposition we find the following expressions for the (almost) Hamiltonian vector fields
associated to the coordinate functions on $\mathcal{R}$ with respect to the standard reduced bracket $\{\cdot,\cdot \}_\mathcal{R}$:
\begin{equation}
\label{E:coord_vector_fields_standard_bracket}
X^\mathcal{R}_{K^b_i}=c_{ijl} \left (  K^b_l +mr^2(  \omega^b_l - \omega^s_3\gamma_l )\right )\partial^\mathcal{R}_{K_j^b}
 + c_{ijl} \gamma_l\partial^\mathcal{R}_{\gamma_j},  \qquad X^\mathcal{R}_{\gamma_i}=-c_{ijl} \gamma_l\partial^\mathcal{R}_{K_j^b},
\end{equation}
where $\partial^\mathcal{R}_{K_j^b}$ and  $\partial^\mathcal{R}_{\gamma_j}$ denote derivations on $\mathcal{R}$ with respect to ${K_j^b}$ and $\gamma_j$.
We are now ready to show

\begin{theorem}
\label{T:Non_integrability_of_distribution}
The characteristic distribution of the reduced standard bracket, $\mathcal{U}:=\{ X_F^\mathcal{R}\; : \; F\in C^\infty(\mathcal{R})\}\subset T\mathcal{R}$, is non-integrable.  
\end{theorem}
\begin{proof}
 At every point in $\mathcal{R}$, any such vector field $X_F^\mathcal{R}$ is a linear combination of
the six vector fields defined in (\ref{E:coord_vector_fields_standard_bracket}). At a generic point in $\mathcal{R}$ only four of them are linearly independent. To see this, first recall that they are all annihilated by
$\frac{1}{2}d||\gamma||^2$ since they are vector fields on $\mathcal{R}$. In addition, a direct calculation shows that they are also annihilated by
the one-form:
\begin{eqnarray*}
\alpha: &=& dK_3^s+mr^2\omega^b\cdot d\gamma=d({\bf K}^b\cdot \gamma) +mr^2\omega^b\cdot d\gamma \\
&=& (K_i^b +mr^2\omega^b_i)d\gamma_i + \gamma_idK^b_i.
\end{eqnarray*}
Thus,  $\alpha$ annihilates any vector in $\mathcal{U}$.
The crucial point is that $\alpha$ is not closed. To formally show that $\mathcal{U}$ is non-integrable we will prove that
\begin{equation}
\label{E:proof_of_nonintegrability}
d\alpha \left ( X^\mathcal{R}_{K^b_i}, X^\mathcal{R}_{\gamma_i} \right )>0.
\end{equation}
Suppose for the moment that the above inequality holds.
In view of the well known identity $d\alpha(X,Y)=X(\alpha(Y))-Y(\alpha(X))-\alpha([X,Y])$, and since $\alpha$ 
annihilates both $X^\mathcal{R}_{K^b_i}$ and $X^\mathcal{R}_{\gamma_i}$, inequality (\ref{E:proof_of_nonintegrability})   implies that there exists $i_0\in\{1,2,3\}$ such that 
\begin{equation*}
\alpha \left ( \left [ X^\mathcal{R}_{K^b_{i_0}}, X^\mathcal{R}_{\gamma_{i_0}} \right ] \right )  \neq 0, \qquad
\mbox{(no sum over $i_0$),}
\end{equation*}
and non-integrability  follows from Frobenius' theorem.

To show that (\ref{E:proof_of_nonintegrability}) holds, notice that in view of (\ref{E:Omega_in_terms_of_K_and_Gamma}) we can write $\omega^b_i=T(\gamma)_{ij}K^b_j$ where
$T(\gamma)_{ij}$ are the components of the $\gamma$ dependent, $3\times 3$ matrix
\begin{equation*}
T(\gamma)=A^{-1}+\left( \frac{mr^2}{1-mr^2(\gamma\cdot A^{-1}\gamma)} \right ) (A^{-1}\gamma)(A^{-1}\gamma)^t.
\end{equation*}
It is clear from the above expression that $T(\gamma)$ is symmetric and positive definite. We can then write
\begin{equation*}
d\alpha=mr^2 d\omega^b_i\wedge \gamma_i=mr^2\left( T(\gamma)_{ij}d K^b_j\wedge d\gamma_i + K^b_j\frac{\partial T(\gamma)_{ij}}{\partial \gamma_k}d\gamma_k\wedge d\gamma_i \right ).
\end{equation*}
Therefore, a direct calculation using (\ref{E:coord_vector_fields_standard_bracket}) gives
\begin{equation*}
d\alpha \left ( X^\mathcal{R}_{K^b_l}, X^\mathcal{R}_{\gamma_l} \right )=T(\gamma)_{ij}c_{lik}c_{ljr}\gamma_k\gamma_r.
\end{equation*}
Using the identity $c_{lik}c_{ljr}=\delta_{ij}\delta_{kr}-\delta_{ir}\delta_{jk}$, and $||\gamma||^2=1$
we get:
\begin{equation*}
d\alpha \left ( X^\mathcal{R}_{K^b_l}, X^\mathcal{R}_{\gamma_l} \right )=\textup{trace}\left( T(\gamma)\right ) - T(\gamma)_{ij}\gamma_i\gamma_j.
\end{equation*}
The above quantity is strictly positive since $T(\gamma)$ is  positive definite  and $||\gamma||=1$.
\end{proof}

As it was mentioned at the end of section \ref{S:Reduction_of_Almost_Poisson_Brackets}, the non-integrability of the characteristic distribution
 implies that there cannot exist a conformal factor for the reduced standard bracket that Hamiltonizes  the problem. As we shall see, the situation is quite different for the reduced affine bracket.

 \subsubsection*{The Reduced Affine Bracket}

 The following proposition gives explicit formulae for the reduced affine nonholonomic bracket $\{\cdot ,\cdot \tilde{\}}_\mathcal{R}$ in  the coordinates $({\bf K}^b,\gamma)$ for the reduced space $\mathcal{R}$.

\begin{proposition} 
\label{P:Formulae_for_reduced_affine_bracket}
We have
\begin{equation*}
\{K^b_i,K^b_j\tilde{\}}_\mathcal{R}=-c_{ijl} \left (  K^b_l - mr^2  \omega^s_3\gamma_l  \right ), \qquad \{K^b_i,\gamma_j \tilde{\}}_\mathcal{R}= -c_{ijl} \gamma_l,  \qquad \{\gamma_i,\gamma_j \tilde{\}}_\mathcal{R}=0.
\end{equation*}
\end{proposition}
The remark made after the statement of proposition \ref{P:Formulae_for_reduced_standard_bracket}
also applies here.
\begin{proof}
The proof is  analogous to that of proposition \ref{P:Formulae_for_reduced_standard_bracket}.
The crucial point is to derive the identities,
 \begin{equation*}
\tilde X_{K^b_i}^\C=\hat X_i^{\textup{left}}+r(g_{2i} \hat \partial_x-g_{1i}\hat \partial_y)+c_{ijl}(K_l^b -mr^2  \omega_3^s \gamma_l ))\hat \partial_{K_j^b}, \qquad \tilde X_{\gamma_i}^\C=c_{ijl}\gamma_l \hat \partial_{K_j^b}.
\end{equation*}
using expression (\ref{E:Omega_Tilde_C_in_body_coordinates}) for $\tilde \Omega_\C$.
\end{proof}

Using this expressions for the bracket we can now show: 
\begin{theorem}
The angular momentum with respect to the vertical axis, $K^s_3={\bf K}^b\cdot \gamma$ is a Casimir function of the affine reduced bracket $\{\cdot,\cdot \tilde{\}}_\mathcal{R}$, i.e. $\{F , {\bf K}^b\cdot \gamma \tilde {\}}_\mathcal{R}=0$ for all $F\in C^\infty(\mathcal{R})$.
\end{theorem}
\begin{proof}
A simple calculation shows $\{K_j^b, {\bf K}^b\cdot \gamma \tilde{\}}_\mathcal{R}=\{\gamma_j, {\bf K}^b\cdot \gamma \tilde{\}}_\mathcal{R}=0$.
\end{proof}
Therefore, in contrast with theorem \ref{T:Non_integrability_of_distribution} we have,
\begin{corollary}
\label{C:Integrability_of_affine_characteristic_distribution}
The characteristic distribution of the reduced  affine bracket, $\mathcal{U}=\{ \tilde X_F^\mathcal{R}\; : \; F\in C^\infty(\mathcal{R})\}$, is everywhere tangent to the foliation 
of $\mathcal{R}$ defined by the level sets of $K_3^s= {\bf K}^b\cdot \gamma$.
\end{corollary}
So the properties of the two reduced brackets are fundamentally different.

\subsection*{Hamiltonization of the Reduced Affine Bracket}

An even stronger result than the one given in corollary \ref{C:Integrability_of_affine_characteristic_distribution} is that the strictly positive function $\mu:Q/G \cong S^2\rightarrow \R$ given by $\mu( \gamma)=Y(\gamma)^{1/2}$, with $Y(\gamma)=1-mr^2(\gamma\cdot A^{-1}\gamma)$, is a conformal factor for the affine reduced bracket.  Define the new bracket $\{\cdot,\cdot \tilde{\}}^\mu_\mathcal{R}$ of functions on $\mathcal{R}$ by the rule:
\begin{equation}
\label{E:def_scaled_bracket}
\{F_1,F_2 \tilde{\}}^\mu_\mathcal{R}:=\mu
\{F_1,F_2 \tilde{\}}_\mathcal{R}.
\end{equation}
This is exactly the bracket for the Chaplygin sphere problem given by the authors in \cite{BorisovMamaev,BorisovMamaev1995}. 
\begin{theorem}
\label{T:Hamiltonization} 
The bracket $\{\cdot,\cdot \tilde{\}}^\mu_\mathcal{R}$ of functions on $\mathcal{R}$ 
defined by (\ref{E:def_scaled_bracket}) satisfies the Jacobi identity.
\end{theorem}

The proof  is a long calculation that will not be included due to space constraints. We can provide the details upon request.
Define a new time $\tau$ by the rescaling:
\begin{equation*}
\label{E:time_change}
dt=\mu \; d\tau.
\end{equation*}
The reduced equations of motion can be written in \emph{Hamiltonian form} in the new time $\tau$
as:
\begin{equation}
\label{E:Ham_form_of_the_reduced_eqns}
\frac{dF}{d\tau}=\{F,\Ham_\mathcal{R} \tilde{\}}^\mu_\mathcal{R}, \qquad \mbox{for all} \qquad F\in C^\infty(\mathcal{R}).
\end{equation}

Moreover, the function  $M^s_3={\bf K}^b\cdot \gamma$ is also a Casimir function for the scaled bracket  $\{\cdot,\cdot \tilde{\}}^\mu_\mathcal{R}$. So the above equation defines a two degree of freedom
Hamiltonian system in each level set of $M^s_3$.
 
\subsection*{The Reduced Equations of Motion and their Integrability}

We now write explicitly the equations of motion and discuss their integrability  in the context of the Hamiltonization discussed above. By differentiating the reduced Hamiltonian 
 (\ref{E:Reduced_Hamiltonian}) one finds after a  lengthy but straightforward calculation,
 \begin{equation*}
\frac{\partial \Ham_\mathcal{R}}{\partial K^b_i}=\omega_i^b, \qquad \frac{\partial \Ham_\mathcal{R}}{\partial \gamma_i}=mr^2\omega_3^s(\omega^b_i-\omega_3^s\gamma_i).
\end{equation*}
Using these expressions and any of the reduced nonholonomic brackets (either the standard or the affine) one computes the reduced equations of motion,
\begin{equation*}
\dot K^b_i=-c_{ijl}K_l^b\omega_j^b, \qquad \dot \gamma_i=c_{jil}\gamma_l\omega_j^b, \qquad \mbox{where}
\qquad \dot{\;}\;=\frac{d}{dt}.
\end{equation*}
In vector form we obtain  the equations that are usually found in the literature:
\begin{equation*}
\label{E:Classical_Chaplygin_Equations}
\dot {\bf K}^b = {\bf K}^b \times \omega^b, \qquad \dot \gamma= \gamma \times \omega^b.
\end{equation*}
These equations have the geometric integral $|| \gamma||^2=1$, the conserved quantity arising from the
$SE(2)$ symmetry, ${\bf K}^b\cdot \gamma$,  and the energy integral, $\Ham_\mathcal{R}$. In addition, the function $J:=\frac{1}{2}{\bf K}^b\cdot {\bf K}^b=\frac{1}{2}\delta_{ij}K_i^bK_j^b$ is directly seen to be in involution with $\Ham_\mathcal{R}$ (with respect to
any of the brackets in $\mathcal{R}$).  In addition to these integrals, one can show that the measure  $\mu(\gamma)^{-1}d{\bf K}^bd\gamma$ is preserved by the flow. It follows that the system is integrable by
quadratures by Jacobi's theorem on the last multiplier. 

The reduced equations were first solved by Chaplygin,
\cite{chapsphere}, in terms of hyper-elliptic functions. A summary of the integrability can be found in \cite{ArnoldIII,FedorovKozlov1995}
were it is shown that the solutions define rectilinear nonuniform motion in two dimensional tori.
The algebraic integrability of the system is considered in \cite{Duistermaat} and a complete complex
solution can be found in \cite{FedorovChaplygin}.

 In (\ref{E:Ham_form_of_the_reduced_eqns}) the equations of motion were written in Hamiltonian form
 in the new time $\tau$. In each symplectic leaf, defined as a  level set of ${\bf K}^b\cdot \gamma$, we have a two-degree of freedom Hamiltonian system. Since $\Ham_\mathcal{R}$ and $J$ are in involution, and their level
sets are compact, we have an integrable Hamiltonian system in the Liouville sense. 

From Liouville's theorem we recover the results given in \cite{ArnoldIII,FedorovKozlov1995} of  uniform rectilinear motion in the new time $\tau$ on two-dimensional tori.  Notice that the existence of a preserved measure follows directly from the Hamiltonization of the problem, it is a multiple of the Liouville measure for the rescaled Hamiltonian flow on each symplectic leaf.  In fact, any (almost) Hamiltonian vector field with respect to the
reduced affine nonholonomic bracket $\{\cdot , \cdot \tilde{\}}_\mathcal{R}$ will preserve the same measure.
In particular, this is the case for the vector field $\tilde X_J^\mathcal{R}$ that by a direct calculation can be shown to be given by
 \begin{equation*}
\tilde X_J^\mathcal{R}=-c_{ijl}K_i^b(A\omega^b)_l-c_{ijl}\gamma_i (A\omega^b)_l.
\end{equation*}
Moreover, since $\Ham_\mathcal{R}$ and $J$ are in involution, the vector fields $\tilde X_J^\mathcal{R}$
and  $\tilde X_{\Ham_{\mathcal{R}}}^\mathcal{R}$ commute after scaling them by $\mu(\gamma)$. The observation that these two vector fields commute already appears in  \cite{Duistermaat} where no reference
to the Hamiltonization of the problem is made.

\section{Final Remarks}
\label{S:Final_Remarks}

To obtain  the Hamiltonization of the Chaplygin sphere problem by reduction  we were forced to introduce the notion of affine  almost Poisson brackets. At this point the presence of the particular affine  form $\Omega_0$ given by (\ref{E:def_Omega_0_for_sphere}) that defines the ``correct" bracket remains a mystery. 
The question remains open to give
a useful characterization of this  form 
in a more general setting.

A possible approach is to consider the 
 affine almost symplectic counterpart. In broad lines, this approach generalizes  the theory of reduction given
in \cite{Koiller1992, BS93,  PlanasBielsa2004} by allowing the formulation to be made in terms of
an affine almost symplectic structure. This approach has the flavor of reduction by stages, it 
sheds some light on the need of the affine term $\Omega_0$, and is part of the content of \cite{NaranjoHochgerner2008}.

Another approach is to consider the reduction of nonholonomic systems using Dirac structures as  was recently
developed in \cite{JotzRatiu2008}. The presence of the affine term  $\Omega_0$ could be related to the theory
of Poisson geometry 
with a 3-form background as introduced in \cite{SeveraWeinstein}. It seems to be a rather strong coincidence that the building block to define the form $\Omega_0$ is precisely the Cartan three-form on $SO(3)$.

We end up by stressing that the key property of the affine Poisson structure that we have considered, is
that the conserved quantity $K_3^s$ becomes a Casimir function of the reduced bracket. This  was
not the case with the standard nonholonomic bracket. It is thus natural to ask the following question: 
 Suppose that a Lie group $H$  is a symmetry group of a nonholonomic system and  that there are conserved quantities associated with its action. Suppose in addition that these conserved quantities are invariant under the action. Does there exist a (possibly affine) nonholonomic bracket for the system such that the conserved quantities are Casimir functions for the corresponding reduced bracket? This issue is also treated in  \cite{NaranjoHochgerner2008}.

\section*{Acknowledgments}

I would like to thank H. Flaschka,  Yu. Fedorov and S. Hochgerner for long and interesting conversations,  and  J. Koiller for 
his encouragement and important remarks reviewing my thesis.

\end{document}